\documentclass[fleqn,12pt,twoside]{article}
\usepackage[headings]{espcrc1}
\usepackage{amsmath}
\usepackage{graphicx}
\PassOptionsToPackage{ctagsplt,righttag}{amsmath}
\usepackage{multicol}
\usepackage{lscape}
\usepackage{listings}

\newcommand{\bc}{\begin{center}}
\newcommand{\ec}{\end{center}}
\newcommand{\beq}{\begin{equation}}
\newcommand{\eeq}{\end{equation}}
\newcommand{\bea}{\begin{eqnarray}}
\newcommand{\eea}{\end{eqnarray}}

\newcommand{\ddx}[2]{\frac{\partial #1}{\partial #2}}
\newcommand{\DDx}[2]{\frac{D #1}{D #2}}

\def\div{\ensuremath{\nabla \cdot}}
\def\divs{\ensuremath{\nabla_{s} \cdot}}
\def\grad{\ensuremath{\nabla}}
\def\grads{\ensuremath{\nabla_{s}}}
\def\lapl{\ensuremath{\nabla^2}}

\def\jac{\ensuremath{{\bf J}}}

\def\diag{\ensuremath{\mathrm{diag}}}
\def\dotstar{\ensuremath{\ .\!*\ }}
\def\dotslash{\ensuremath{\ ./\ }}
\def\dotwedge{\ensuremath{.\!\wedge}}

\def\eg{\emph{e.g., }}
\def\ie{\emph{i.e., }}

\def\ch{\ensuremath{\hat{c}}}
\def\psih{\ensuremath{\hat{\psi}}}

\def\qh{\ensuremath{\hat{q}}}
\def\wh{\ensuremath{\hat{w}}}

\begin{document}
\bibliographystyle{elsart-num-sort}

\title{A Direct Matrix Method for Computing Analytical Jacobians of 
       Discretized Nonlinear Integro-differential Equations}

\author{
Kevin T. Chu\address{Vitamin D, Inc., Menlo Park, CA 94025}$^,$\address{Institute of High Performance Computing, A*STAR, Singapore, Singapore}
}

\runtitle{Direct Matrix Method for Computing Jacobians of Discretized Integro-diff. Eqs.}
\runauthor{K. T. Chu}

\maketitle

\noindent \rule{6.3in}{1pt}

\begin{abstract}
In this pedagogical article, we present a simple direct matrix method for 
analytically computing the Jacobian of nonlinear algebraic equations that 
arise from the discretization of nonlinear integro-differential equations.  
The method is based on a formulation of the discretized equations in vector 
form using only matrix-vector products and component-wise operations.  By 
applying simple matrix-based differentiation rules, the matrix form of the 
analytical Jacobian can be calculated with little more difficulty than that 
required when computing derivatives in single-variable calculus.
After describing the direct matrix method, we present numerical experiments 
demonstrating the computational performance of the method, discuss its 
connection to the Newton-Kantorovich method, and apply it to illustrative 
1D and 2D example problems.  MATLAB code is provided to demonstrate the 
low code complexity required by the method.
\end{abstract} \newline 

\noindent {\footnotesize \emph{Keywords:} analytical Jacobian; 
numerical methods; matrix calculus; Newton's method; 
integro-differential equations}

\noindent \rule{6.3in}{1pt}


\section{Introduction} 
Many numerical methods for solving nonlinear integro-differential equations 
require computation of the Jacobian for the system of algebraic equations 
that arises when the continuous problem is discretized.  For example, 
any Newton's method calculation requires computation of the Jacobian (exactly 
or approximately) during each Newton 
iteration~\cite{kelley_iterative_methods,kelley_newtons_method}.  
Unfortunately, calculation of the Jacobian can be a time-consuming and 
error-prone procedure for both the computer \emph{and} the scientific 
programmer.  

In this pedagogical article, we present a simple \emph{direct matrix method} 
for calculating analytical Jacobians of discretized, nonlinear 
integro-differential equations.  The direct matrix method produces the 
Jacobian for the discretized equations directly in matrix form \emph{without} 
requiring calculation of individual matrix elements.  The essential idea is
to first write the discretized, integro-differential equation explicitly in 
terms of discrete operators (\eg differentiation and quadrature 
matrices~\cite{boyd_spectral_book,fornberg_spectral_book,trefethen_spectral_book}) 
and then use simple \emph{matrix}-based differentiation rules to calculate 
the Jacobian directly~\cite{chen1999,chen2000a,chen1997a,chu_thesis_2005}. 
The key observation underlying this approach is that there is 
a tremendous amount of structure in the nonlinear algebraic equations that 
arise from the discretization of nonlinear integro-differential equations.  
By taking advantage of this structure, the calculation of analytical Jacobians
is reduced to nearly the level of complexity required to compute derivatives 
of scalar, single variable functions.

The operator-based approach of expressing and analyzing discretized differential
equations has been used implicitly by the scientific computing community for 
decades, especially in the context of the Newton-Kantorovich and related
methods~\cite{boyd1986b,boyd1986a,kelley1987}.  However, a direct matrix
approach seems to have been first formally described in one place by Chen who 
presented a collection of rules formulated in terms of specially defined 
matrix products~\cite{chen1999,chen2000a,chen1997a}.  
In addition to using the method to solve nonlinear partial differential 
equations~\cite{chen2000a,chen1998}, Chen used his formulation of 
the method to develop several interesting theoretical results (\eg
stability analysis of numerical methods for nonlinear time-dependent problems) 
based on the observation that when a nonlinear differential equation only has 
polynomial nonlinearities, there is a very close relationship between the 
discretized nonlinear differential equation and its 
Jacobian~\cite{chen1999,chen2000b}.  

Mathematically, the present formulation of the direct matrix method is 
equivalent to Chen's approach.  However, rather than introducing special
matrix products, we rely solely on standard linear algebra operations 
augmented by component-wise operations (\eg the Hadamard or Schur 
product~\cite{horn_topics}).  In addition, we have chosen to use MATLAB 
notation in our formulation because of its prevalence in modern scientific 
computing.  Working in MATLAB notation has the added benefit of making it 
almost trivial to translate the analytical calculations into working numerical 
code\footnote{With today's powerful desktop and laptop computers, MATLAB
is quite capable of handling moderate-sized production work.}.

Another feature of our formulation, which is also present to some extent 
in~\cite{chen1999}, is the emphasis on the analogy between calculation of 
Jacobians for discretized, nonlinear integro-differential equations and 
calculation of derivatives for scalar functions of a single variable.  To 
help strengthen the analogy with single-variable calculus, we organize the 
operations required to compute a Jacobian as a short list of simple 
differentiation rules.  

This article is organized as follows.  In the remainder of this section,
we compare the direct matrix method with several common methods for computing 
Jacobians.  In Section~\ref{sec:direct_matrix_method}, we present the
direct matrix method, including a discussion of its computational performance
and its relation to the Newton-Kantorovich method.  
Finally, in Section~\ref{sec:applications}, we apply the direct matrix method 
to two examples (one 1D and one 2D) from the field of electrochemical 
transport.  To demonstrate the low code complexity required by the direct 
matrix method, MATLAB code for the example problems is provided in the 
appendices. 
Throughout our discussion, we will focus solely on collocation methods where 
the continuous and discrete forms of the integro-differential equation have 
essentially the same structure.  However, it is important to recognize that 
the direct matrix method can also be used for Galerkin methods by applying 
it directly to the weak-form of the problem.

\subsection{Comparison with Common Methods for Computing Jacobians}
One common approach for obtaining the Jacobian of a discretized, nonlinear 
integro-differential equation is to compute it numerically using finite 
differences of the grid function or expansion coefficient 
values~\cite{kelley_iterative_methods,kelley_newtons_method}. 
Unfortunately, numerical computation of the Jacobian can be time consuming 
for some problems.  Depending on the numerical method, it might be possible 
to reduce the computational cost of a numerical Jacobian by taking advantage 
of the sparsity pattern in the 
Jacobian~\cite{coleman1983,curtis1974,kelley_newtons_method}, but determining
the sparsity pattern can be complicated for nontrivial problems.

As an example, consider the Poisson-Nernst-Planck equations for 
electrochemical transport~\cite{bazant2005,chu2006,rubinstein_book}:
\bea 
\ddx{c_+}{t} &=& \div \left( \grad c_+ + c_+ \grad \phi \right) 
  \label{eq:cation_NP} \\
\ddx{c_-}{t} &=& \div \left( \grad c_- - c_- \grad \phi \right) 
  \label{eq:anion_NP} \\
\epsilon \lapl \phi &=& - \left( c_+ - c_- \right)
  \label{eq:poisson_eqn}, 
\eea 
where $c_\pm$ are cation and anion concentrations, respectively, and
$\phi$ is the electric potential, and $\epsilon$ is a dimensionless physical 
parameter related to the dielectric constant of the electrolyte.  
The first two of these are the Nernst-Planck equations for ion transport 
and are simply conservation laws for cations and 
anions~\cite{newman_book,rubinstein_book}.  
The last equation is the Poisson equation~\cite{jackson_book}, which provides
closure for the Nernst-Planck equations.  Note that in (\ref{eq:poisson_eqn}),
the local charge density has been written in terms of the individual
ions, which are the only source of charge density in many electrochemical 
systems.
(\ref{eq:cation_NP}) -- (\ref{eq:poisson_eqn}) form a nonlinear parabolic
system of partial differential equations, which suggests that we use
an implicit time-stepping scheme to numerically solve the equations.
This choice, however, requires that at time $t_{n+1}$, we solve a 
nonlinear system of equations for $c_\pm^{(n+1)}$ which depends on an 
auxiliary variable $\phi$ which is in turn related to $c_\pm^{(n+1)}$
through the Poisson equation.  Note that in order to numerically compute 
the Jacobian for the resulting nonlinear system of equations for 
$c_\pm^{(n+1)}$, we must solve the Poisson equation for each perturbation 
to the current iterate of $c_\pm^{(n+1)}$.  Therefore, for a pseudospectral 
discretization of (\ref{eq:cation_NP}) -- (\ref{eq:poisson_eqn}) using $N$ grid 
points, numerically computing the Jacobian requires $O(N^4)$ operations, which 
is much higher than the $O(N^2)$ elements in the Jacobian\footnote{While 
low-order discretization of the equations do not show this same disparity in 
the computation time and the number of elements (requiring $O(N^2)$ operations 
for $O(N^2)$ elements), they typically require many more grid points to produce
an accurate solution.}.

Using an analytical Jacobian is one way to avoid the computational cost 
associated with numerical Jacobians.  In principle, it is straightforward to 
derive the analytical Jacobian for the system of algebraic equations that 
arises when a nonlinear integro-differential equation is discretized.  
Index notation (also known as tensor notation) is perhaps the most common 
technique used to calculate analytical Jacobians.  The basic idea behind the 
index notation method is to write the discretized form of the differential 
equation using index notation and then use tensor calculus to compute 
individual matrix elements in the Jacobian.  For example, for a 
finite-difference or pseudospectral discretization, the discretized equations 
can be written in the form:
\beq
  F_i(u_1, u_2, \ldots, u_N) = 0 
  \label{eq:index_notation_eqns}
\eeq
for $i = 1, 2, \ldots ,N$ where $u_i$ and $F_i$ are the value of the solution
and the discretized differential equation at the $i$-th grid point (or more 
generally, the $i$-th collocation point).  
Boundary conditions are included in this formulation by using the discretized 
boundary conditions (rather than the integro-differential equation) at grid 
points on the boundary (or immediately adjacent to if no grid points reside
on the boundary)\footnote{Care must be exercised when imposing boundary conditions, especially when using pseudospectral methods~\cite{fornberg_2006}.}.  
The $ij$-th element of the Jacobian, $\jac$, for 
(\ref{eq:index_notation_eqns}) is simply the partial derivative of $F_i$ 
with respect to $u_j$: $\jac_{ij} = \partial F_i / \partial u_j$.  
While simple and straightforward, index notation suffers from the 
disadvantage of being somewhat tedious and error-prone.  The main
challenge in using index notation is keeping track of all of the indices
when writing out and computing partial derivatives of the discrete
equations.  

Automatic differentiation~\cite{griewank_book,autodiff_anl} offers an 
important alternative when exact Jacobians are desired.  Because it generates 
code for computing the Jacobian directly from the code used to evaluate the 
residual, automatic differentiation completely eliminates the possibility of 
human error when deriving the exact Jacobian and implementing it in code.  
Recent developments have made automatic differentiation available in several 
common programming languages (including MATLAB~\cite{shampine2005}).  While 
useful, automatic differentiation still takes some effort to use and may not 
always generate the most compact, efficient code.  Active development in this 
area will certainly continue to improve the usability of automatic 
differentiation software and the performance of generated code. 

The direct matrix method has several advantages over the methods discussed
in this section.  First, the direct matrix method yields a more accurate
Jacobian than finite differences and generally in less time (see 
Section~\ref{sec:computational_perf}). 
Second, because the method is based on simple differentiation rules, the 
calculation is straightforward and less prone to error than the index notation 
approach.  The differentiation rules also make it easier to calculate the 
Jacobian for differential equations which depend on auxiliary variables, such 
as (\ref{eq:cation_NP}) -- (\ref{eq:poisson_eqn}).
From a programming perspective, calculation of the Jacobian directly in matrix 
form facilitates implementation of numerical methods for nonlinear problems in 
languages that have built-in support for matrix and vector operations 
(\eg MATLAB and Fortran 95).  Finally, having the Jacobian available in matrix 
form can be useful for analyzing properties of numerical 
methods~\cite{chen1999}.  While it may be possible to to convert the 
element-wise representation of the Jacobian derived using index notation 
for simple problems, this step can be challenging for more complex 
problems\footnote{Interestingly, the conversion from element-wise to matrix 
representation of the Jacobian often reveals the close relationship between 
the Jacobian for the discrete equations and the underlying structure of the 
original integro-differential equation.}.

\section{The Direct Matrix Method\label{sec:direct_matrix_method}}
There are two basic ideas underlying the direct matrix method for 
calculating analytical Jacobians of discretized, integro-differential 
equations.  First, rather than writing the discretized, integro-differential 
equations at each of the collocation points in terms of individual elements 
of the solution vector, we write the entire system of equations as a single 
vector equation expressed explicitly in terms of matrix-vector products and 
component-wise multiplication (\eg Hadamard products).  Second, the analytical 
Jacobian for the discretized system of equations is computed directly in 
matrix form by using simple differentiation rules that are reminiscent of 
those used to compute derivatives in single-variable calculus.  In this 
section, we develop both of these ideas in detail.  Towards the end of the 
section, we comment on the computational performance of the direct matrix
method and its relationship to the Newton-Kantorovich 
method~\cite{boyd_spectral_book}.

\subsection{Matrix-vector Representation of Discretized Equations}
Writing the discretized, nonlinear integro-differential equation explicitly 
in terms of basic linear algebra and component-wise algebraic operations 
is the initial step of the direct matrix method.  Because of the
similarities in the structure between the discrete and continuous forms
of the equations, the procedure is very straightforward.
First, convert all differentiation and integration operators into their 
discrete analogues.  Since both of these operations are linear, they 
become multiplication of vectors representing grid functions by 
differentiation and quadrature matrices, respectively:
\bea
 \frac{du}{dx} &\rightarrow&  D * \hat{u} \\
 \int u dx &\rightarrow&  Q * \hat{u},
\eea
where the hat accent indicates a discretized field variable
and $D$ and $Q$ are the differentiation and quadrature matrices
associated with the choice of computational grid. 

Next, convert all point-wise algebraic operations and function evaluations
in the continuous equations to component-wise algebraic operations and
function evaluations in the discrete equations.
Some examples of the conversion process include:
\bea 
  u \frac{dv}{dx} &\rightarrow& 
    \hat{u} \dotstar \left( D*\hat{v} \right) \\
  \sin(u) &\rightarrow& \sin(\hat{u}) \\
  u^2 &\rightarrow& \hat{u} \dotwedge 2 \\
  \exp(u) &\rightarrow& \exp(\hat{u}).
\eea 
In these examples, we have adopted the MATLAB convention of using $. op$ to 
represent component-wise application of the $op$ operation.  Also, note that 
we have abused notation for component-wise function evaluations --
$f(\hat{u})$ represents the vector 
\beq
(f(\hat{u}_1), f(\hat{u}_2), \ldots, f(\hat{u}_N))
\eeq 
not an arbitrary vector function of the entire solution vector $\hat{u}$.
Throughout our discussion, we will indicate component-wise and general
functions of $\hat{u}$ by using lowercase and uppercase variables, 
respectively.

\subsubsection{Differential and Integral Operators in Multiple Space 
               Dimensions \label{sec:multiple_space_dims}}
It is important to emphasize that the matrix-vector representation is
not restricted to scalar field equations or functions of a single 
variable.  Handling vector equations is simple -- vector equations 
may be treated as systems of equations and vector operations may be 
expressed in component-wise form.  

The construction of differential and integral operators for functions of 
multiple variables is slightly more complicated but straightforward.  First, 
we represent grid functions as a 1D vector by selecting an ordering of the grid 
points.  Then, we derive the differentiation and integration matrices 
associated with this choice of ordering.  While there is no unique mapping 
from a multi-dimensional grid to a 1D vector, it is important to choose the 
ordering of the grid points carefully because it directly affects the ease 
with which differentiation and integration matrices can be derived.

\begin{figure}
\bc
\scalebox{0.75}{\includegraphics{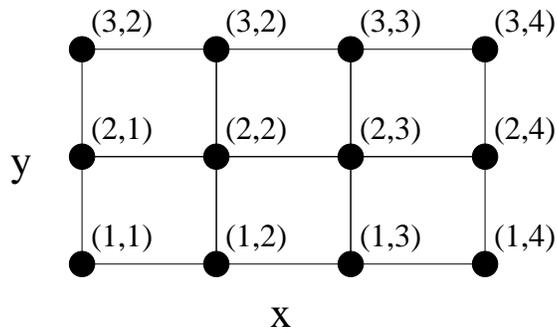}}
\caption[Example of a 2d grid]{
\label{figure:2d_example_grid}
Example of a 2d grid with four grid points in the $x$-direction and
three grid points in the $y$-direction.
}
\ec
\end{figure}

For problems on logically rectangular computational domains, the computational 
grid may be constructed as a Cartesian product of one-dimensional 
grids~\cite{trefethen_spectral_book}.  The 
result is a structured grid (possibly non-uniform depending on the 
discretization in each coordinate direction).  To represent a grid function 
as a 1D vector, the most natural way to flatten the grid is by using a 
lexicographic order for the grid indices.  For example, on the small 
$3$ by $4$ Cartesian grid in Figure~\ref{figure:2d_example_grid}, we could 
order the grid function, $\hat{u}$, one row at a time:
\beq
  \hat{u} = \left(
    \hat{u}_{11}, \hat{u}_{12}, \hat{u}_{13}, \hat{u}_{14},
    \hat{u}_{21}, \hat{u}_{22}, \hat{u}_{23}, \hat{u}_{24},
    \hat{u}_{31}, \hat{u}_{32}, \hat{u}_{33}, \hat{u}_{34}
  \right)^T
\eeq
where $\hat{u}_{ij}$ is the value of $u$ at $(x_i,y_j)$.
With this choice of ordering, discrete partial derivative operators
are given as Kronecker product of the differentiation matrices with 
identity matrices~\cite{trefethen_spectral_book}.
For the example in Figure \ref{figure:2d_example_grid}, the
\emph{partial} differentiation matrices are
\bea
  D_x &=& I_3 \otimes D_4
  \\
  D_y &=& D_3 \otimes I_4
\eea
where $D_n$ and $I_n$ are the 1d differentiation and identity matrices of 
size $n$ and $\otimes$ denotes the Kronecker product.  Note that care must be 
taken to ensure that the order of the Kronecker products is consistent with 
the ordering of the grid function vector.  

To illustrate these ideas, let us consider the discretized forms of the curl 
and divergence of a vector field, $\vec{U}$.  Once we have chosen an ordering 
of the grid points and derived the corresponding partial differentiation 
matrices $D_x$, $D_y$, and $D_z$, the discrete curl operator may be easily 
written as
\bea
  \left( 
  \begin{array}{c}
    D_y * U_z - D_z * U_y \\
    D_z * U_x - D_x * U_z \\
    D_x * U_y - D_y * U_x 
  \end{array}
  \right),
\eea
where $U_x$, $U_y$, and $U_z$ are the components of $\vec{U}$.  
The discrete divergence operator is also straightforward to derive:
\bea
D_x * U_x + D_y * U_y + D_z * U_z.
\eea

\paragraph{Boundary Conditions} 
For problems in multiple space dimensions, it is often convenient to include
discretized boundary conditions in a matrix-vector representation by
first breaking the full differentiation and integration matrices into 
multiple components.  Each component is defined by the grid points 
that it contributes to and the grid points it receives contributions from. 
For instance, it may be convenient to decompose a differentiation matrix into 
four mutually exclusive components that: 
(1) use interior points and contribute to interior points, 
(2) use interior points and contribute to boundary points, 
(3) use boundary points and contribute to interior points, 
and 
(4) use boundary points and contribute to boundary points.
Splitting the differentiation and integration matrices into separate 
components can be helpful when computing residuals and Jacobians for the 
different types of grid points in the computational domain.

Deriving these components is straightforward using zeroth-order restriction 
and prolongation matrices~\cite{briggs_book}. 
A zeroth-order restriction matrix is a matrix of zeros and ones which extracts 
a desired subset of elements from a vector.  A zeroth-order prolongation 
matrix is also a matrix of zeros and ones, but it injects the elements of a 
restricted vector into a desired subset of the elements of a full-length 
vector.  If $\{i_1, i_2, \ldots, i_m\}$ are the (flattened) indices of a subset 
of points from $N$ grid points, then the associated restriction matrix, $R$, 
would be an $m \times N$ matrix, $R$, with ones at the positions 
$(1, i_1), (2, i_2), \ldots, (m, i_m)$.  The associated prolongation matrix, 
$P$, that injects a vector of length $m$ into the positions 
$\{i_1, i_2, \ldots, i_m\}$ of a vector of length $N$ 
is simply the transpose of $R$: $P = R^T$.  For example, the restriction
and projection matrices for the two interior points of the grid in 
Figure \ref{figure:2d_example_grid} are given by
\beq
  R = \left[
      \begin{array}{cccccccccccc}
      0 & 0 & 0 & 0 & 0 & 1 & 0 & 0 & 0 & 0 & 0 & 0 \\
      0 & 0 & 0 & 0 & 0 & 0 & 1 & 0 & 0 & 0 & 0 & 0
      \end{array}
      \right]
      \ \ , \ \ P = R^T
\eeq

To derive the differentiation matrix that uses values from grid points
$A$ and contributes to grid points $B$, we use the prolongation matrix 
for set $A$ to expand the restricted grid function associated with $A$ into 
a full-length vector and the restriction matrix for set $B$ to compress
derivative grid function to the set $B$:
\beq
  D_{A\rightarrow B} = R_B * D * P_A = R_B * D * R_A^T.
  \label{eq:restricted_diff_matrix}
\eeq
Restricted integration matrices are derived in exactly the same manner.

\subsection{Simple Differentiation Rules for Computing Exact Jacobians
            \label{sec:differentiation_rules}}
Once the continuous equations have been put in a discretized form that 
is expressed explicitly in terms of matrix-vector products and component-wise 
operations, the analytical Jacobian for the discretized equations 
can be calculated by applying a few simple \emph{matrix}-based
differentiation rules.  Because the differentiation rules are expressed
completely in matrix form without any reference to individual elements
in the Jacobian matrix, they allow us to compute the Jacobian 
directly in matrix form.  
In this section, we list these differentiation rules, which are essentially 
results from multivariate and matrix calculus applied specifically to the 
structure of discretized integro-differential equations.

\subsubsection{Matrix-Vector Product Rule}
The Jacobian of a matrix-vector product (which corresponds to a linear 
operator acting on a function in the continuous equations) is just the 
matrix itself:
\beq
  \frac{\partial}{\partial \hat{u}} \left( A * \hat{u} \right) = A.
  \label{eq:jacobian_mat_vec_rule}
\eeq
For example, the Jacobian of the discretized derivative of $u$, 
$D * \hat{u}$, is just $D$.

\subsubsection{Diagonal Rule}
The Jacobian of a component-wise function $f$ of a grid function $\hat{u}$ is 
a diagonal matrix with diagonal entries given by $f'(\hat{u})$: 
\beq
  \frac{\partial f(\hat{u})}{\partial \hat{u}} =
    \diag \left( f'(\hat{u}) \right)
  \label{eq:jacobian_diag_rule}.
\eeq
In essence, the diagonal rule is a way to use matrix notation to
represent the fact that the differential in the $i$-th component of 
$f(\hat{u})$ only depends on the change in the $i$-th component of $\hat{u}$
and is given by $\delta f(\hat{u}_i) = f'(\hat{u}_i) \delta \hat{u}_i$,
As an example, the Jacobian of $\sin(\hat{u})$ is 
$\diag \left( \cos(\hat{u}) \right)$.

\subsubsection{Chain Rules}
The Jacobian of a matrix $A$ times an arbitrary function, $F$, of \emph{all}
of the components of $\hat{u}$ is $A$ times
the Jacobian of $F$:
\beq
  \frac{\partial}{\partial \hat{u}} \left( A * F(\hat{u}) \right)
    = A * \frac{\partial F}{\partial \hat{u}}
  \label{eq:general_jacobian_chain_rule_1}
\eeq
Similarly, the Jacobian of a function, $F(\hat{u})$, when its argument is 
a matrix $A$ times the grid function $\hat{u}$ is the Jacobian of $F$ 
evaluated at $A*\hat{u}$ times $A$:
\beq
  \frac{\partial}{\partial \hat{u}} F \left( A * \hat{u} \right)
    = \left[\frac{\partial F}{\partial \hat{u}}
      \left(A*\hat{u}\right) \right] * A
  \label{eq:general_jacobian_chain_rule_2}
\eeq
These rules are simply the chain rules for vector fields from 
multivariate calculus~\cite{apostol_book}.

For the special but common case when $F$ is a component-wise function 
$F(\hat{u}) = f(\hat{u})$, (\ref{eq:general_jacobian_chain_rule_1}) reduces 
to $A$ times the diagonal matrix with $f'(\hat{u})$ on the diagonal:
\beq
  \frac{\partial}{\partial \hat{u}} \left( A * f(\hat{u}) \right)
    = A * \diag \left( f'(\hat{u}) \right)
  \label{eq:jacobian_chain_rule_1}
\eeq 
and (\ref{eq:general_jacobian_chain_rule_2}) reduces to the diagonal matrix 
with $f'(A*\hat{u})$ on the diagonal times the matrix $A$:
\beq
  \frac{\partial}{\partial \hat{u}} f \left( A * \hat{u} \right)
    = \diag \left( f'(A*\hat{u}) \right) * A 
  \label{eq:jacobian_chain_rule_2}
\eeq

\subsubsection{Product Rule}
To compute the Jacobian of a component-wise product of general functions 
$F$ and $G$ of a grid function $\hat{u}$, we use the product rule:
\beq
  \frac{\partial}{\partial \hat{u}} 
    \left( F(\hat{u}) \dotstar G(\hat{u}) \right)
    = \diag \left( G(\hat{u}) \right) * \frac{\partial F}{\partial \hat{u}} 
    + \diag \left( F(\hat{u}) \right) * \frac{\partial G}{\partial \hat{u}} 
  \label{eq:jacobian_prod_rule}.
\eeq
The derivation of the product rule follows from the expression for the
variation of the $i$-th component of $F(\hat{u}) G (\hat{u})$
\beq
  \delta \left( F_i(\hat{u}) G_i(\hat{u}) \right) = 
    G_i(\hat{u}) \frac{\partial F_i}{\partial \hat{u}} \delta \hat{u}
  + F_i(\hat{u}) \frac{\partial G_i}{\partial \hat{u}} \delta \hat{u}
  = \left(
    G_i(\hat{u}) \frac{\partial F_i}{\partial \hat{u}} 
  + F_i(\hat{u}) \frac{\partial G_i}{\partial \hat{u}}
    \right) \delta \hat{u},
\eeq
which yields the Jacobian 
\beq
  \left[
  \begin{array}{c}
    G_1(\hat{u}) \left( \partial F_1/\partial \hat{u} \right) \\
    G_2(\hat{u}) \left( \partial F_2/\partial \hat{u} \right) \\
    \vdots \\
    G_N(\hat{u}) \left( \partial F_N/\partial \hat{u} \right) 
  \end{array}
  \right]
  +
  \left[
  \begin{array}{c}
    F_1(\hat{u}) \left( \partial G_1/\partial \hat{u} \right) \\
    F_2(\hat{u}) \left( \partial G_2/\partial \hat{u} \right) \\
    \vdots \\
    F_N(\hat{u}) \left( \partial G_N/\partial \hat{u} \right) 
  \end{array}
  \right]
  = \diag \left( G(\hat{u}) \right) * \frac{\partial F}{\partial \hat{u}}
  + \diag \left( F(\hat{u}) \right) * \frac{\partial G}{\partial \hat{u}}
\eeq

\subsection{Example Jacobian Calculations}
As our first example, let us consider the 1D Poisson equation:
\beq
  \frac{d^2 u}{dx^2} + \rho = 0.
\eeq
To put this in discretized form, we need only replace the continuous second
derivative operator by a discrete analogue:
\beq
  D^2 * \hat{u} + \hat{\rho} = 0.
\eeq
Here, we have chosen to apply the discrete single derivative operator
twice.  Via a direct application of the matrix-vector product rule 
(\ref{eq:jacobian_mat_vec_rule}), the Jacobian of the left-hand side of 
this equation is easily found to be $D^2$.  
Since this is a linear equation, there would not normally be a need to 
compute the Jacobian of the left hand side of this equation.  Moreover, 
the Jacobian for this example is very easy to calculate using alternative 
means (or even by inspection).  We merely present it to illustrate 
the direct matrix method on a simple model problem.

As a less trivial, let us calculate the Jacobian for the discretized 
form of the nonlinear function 
\beq
f(u) = e^{2u} \frac{du}{dx}. 
\eeq
Converting this function to discrete form, we obtain
\beq
f(\hat{u}) = \left( e \dotwedge (2\hat{u}) \right) \dotstar (D * \hat{u}).
\eeq
Using the product rule (\ref{eq:jacobian_prod_rule}), we find that the
the Jacobian is given by 
\beq
  J = \diag \left( D * \hat{u} \right) * 
    \frac{\partial}{\partial \hat{u}} \left(e \dotwedge (2 \hat{u}) \right)
  + \diag \left( e \dotwedge (2 \hat{u}) \right) * 
    \frac{\partial (D * \hat{u})}{\partial \hat{u}}.
\eeq
Then applying the diagonal rule (\ref{eq:jacobian_diag_rule}) and the 
matrix-vector product rule (\ref{eq:jacobian_mat_vec_rule}), we find that
\beq
  J = 2 \diag \left( D * \hat{u} \right) * 
      \diag \left( e \dotwedge (2 \hat{u}) \right) 
    + \diag \left( e \dotwedge (2 \hat{u}) \right) * D,
\eeq
which can be simplified to 
\beq
  J = 2 \diag \left( (D * \hat{u}) \dotstar (e \dotwedge (2 \hat{u})) \right) 
    + \diag \left( e \dotwedge (2 \hat{u}) \right) * D
\eeq
by observing that 
$  \diag \left( \hat{u} \dotstar \hat{v} \right) = 
  \diag \left( \hat{u} \right) * \diag \left( \hat{v} \right). 
$

As a final example, let us calculate the Jacobian for the nonlinear
algebraic equations that arise when solving the one-dimensional version
of (\ref{eq:cation_NP}) -- (\ref{eq:poisson_eqn}) using a simple 
backwards Euler discretization in time.  Using the direct matrix approach
for the spatial discretization, the nonlinear algebraic equations for 
$\hat{c}_+^{(n+1)}$ and $\hat{c}_-^{(n+1)}$ that need to be solved at each 
time step are:
\bea
  \hat{c}_+^{(n+1)} - \Delta t 
  \left (D^2 * \hat{c}_+^{(n+1)} + 
         D * \left(\hat{c}_+^{(n+1)} \dotstar (D*\hat{\phi}) \right) \right)
  -\hat{c}_+^{(n)} &=& 0 \label{eq:discrete_c_plus_eqn} \\
  \hat{c}_-^{(n+1)} - \Delta t 
  \left (D^2 * \hat{c}_-^{(n+1)} - 
         D * \left(\hat{c}_-^{(n+1)} \dotstar (D*\hat{\phi}) \right) \right)
  -\hat{c}_-^{(n)} &=& 0 \label{eq:discrete_c_minus_eqn} \\
  \epsilon D^2*\hat{\phi} + \left(\hat{c}_+^{(n+1)}-\hat{c}_-^{(n+1)}\right) 
    &=& 0 \label{eq:discrete_poisson_eqn} 
\eea
where $\hat{c}_\pm^{(n)}$ are the concentrations at the current time step
and $\Delta t$ is the time step size.  It is important to mention that 
several of the rows in (\ref{eq:discrete_poisson_eqn}) will typically be 
replaced to impose the discretized form of the boundary conditions for 
$\phi$.  For illustrative purposes, let us suppose that we have simple 
Dirichlet boundary conditions for $\phi$.  In this situation, 
(\ref{eq:discrete_poisson_eqn}) is only imposed at interior grid 
points~\cite{trefethen_spectral_book}.  

Using the simple differentiation rules from the previous section, the 
Jacobian of (\ref{eq:discrete_c_plus_eqn}) with respect to $\hat{c}_+^{(n+1)}$ 
is
\beq
  I - \Delta t \left (D^2 
                     + D*\diag\left(D*\hat{\phi}\right) 
                     + D*\diag\left(\hat{c}_+^{(n+1)}\right)*D* 
           \frac{\partial \hat{\phi} \ \ \ \ \ }{\partial \hat{c}_+^{(n+1)}}
  \right),
  \label{eq:J_c_plus}
\eeq
where $I$ is the identity matrix and 
$\left(\partial \hat{\phi} / \partial \hat{c}_+^{(n+1)}\right)$
is the Jacobian of $\hat{\phi}$ with respect to $\hat{c}_+^{(n+1)}$.  To
eliminate
$\left(\partial \hat{\phi} / \partial \hat{c}_+^{(n+1)}\right)$ from
this expression, we simply apply the differentiation rules to 
(\ref{eq:discrete_poisson_eqn}) with two rows eliminated for the 
boundary conditions and solve for the interior portion of
$\left(\partial \hat{\phi} / \partial \hat{c}_+^{(n+1)}\right)$: 
\beq
\left(\frac{\partial \hat{\phi} \ \ \ \ \ }{\partial \hat{c}_+^{(n+1)}}\right)_{int} = - \frac{1}{\epsilon} \left(D^2 \right)_{int}^{-1} 
  \label{eq:d_phi_d_c_plus_int},
\eeq
where $\left(D^2 \right)_{int}$ is the submatrix of $D^2$ that remains when 
all of the columns and rows corresponding to boundary grid points have been 
removed.  
Since the boundary values of $\hat{\phi}$ are fixed and the values of 
$\hat{c}_+^{(n+1)}$ at the boundaries do not affect the potential in the 
interior, the full Jacobian 
$\left(\partial \hat{\phi} / \partial \hat{c}_+^{(n+1)}\right)$ is
given by 
\beq
\frac{\partial \hat{\phi} \ \ \ \ \ }{\partial \hat{c}_+^{(n+1)}} = 
  \left[
  \begin{array}{ccc}
    0 & \cdots & 0 \\
    \vdots & - \frac{1}{\epsilon} \left(D^2 \right)_{int}^{-1} & \vdots \\
    0 & \cdots & 0 \\
  \end{array}
  \right]
  \label{eq:d_phi_d_c_plus},
\eeq
where we have assumed that the first and last grid points correspond to
boundary points.  It is important to recognize that the form for
$\left(\partial \hat{\phi} / \partial \hat{c}_+^{(n+1)}\right)$ 
in (\ref{eq:d_phi_d_c_plus}) is specific to problems with Dirichlet boundary 
conditions for $\phi$.  
For other boundary conditions, the inversion of the equation for 
$\left(\partial \hat{\phi} / \partial \hat{c}_+^{(n+1)}\right)$ generally
leads to different forms for the Jacobian.

The Jacobian for (\ref{eq:discrete_c_plus_eqn}) can now be explicitly 
computed by substituting (\ref{eq:d_phi_d_c_plus}) into (\ref{eq:J_c_plus}).
The similar expression for the Jacobian of (\ref{eq:discrete_c_minus_eqn}) is 
obtained using an analogous procedure.  Using the direct matrix approach, we 
have reduced the calculation of the Jacobian to $O(N^3)$ (cost of 
matrix-inversion and matrix-matrix multiplies) compared to the 
$O(N^4)$ cost for computing a numerical Jacobian for high-order spatial 
discretizations.
It is worth pointing out that in this example, the Jacobian for the 
concentrations does \emph{not} depend explicitly on $\phi$ because the 
Poisson equation is linear.  As a result, there is no need to solve for 
$\phi$ in order to compute the Jacobians for (\ref{eq:discrete_c_plus_eqn}) 
and (\ref{eq:discrete_c_minus_eqn}).  For general problems, the Jacobian
may depend on the auxiliary variable, so it might be necessary to solve the 
constraint equation.  However, because only one solve for the auxiliary 
variables is required with the direct matrix method, the cost of computing 
the Jacobian is still dramatically reduced compared to using finite 
differences.

\subsection{Computational Performance \label{sec:computational_perf}}
In general, using the direct matrix method to compute a Jacobian is faster 
than calculating a numerical Jacobian.  As mentioned in the previous section,
the performance difference is expected to be large when auxiliary variables
are involved in the expression of the residual.  However, the direct matrix
method yields higher performance even for problems where the residual is 
relatively simple.  
\begin{figure}[tb]
\bc
\scalebox{0.35}{\includegraphics{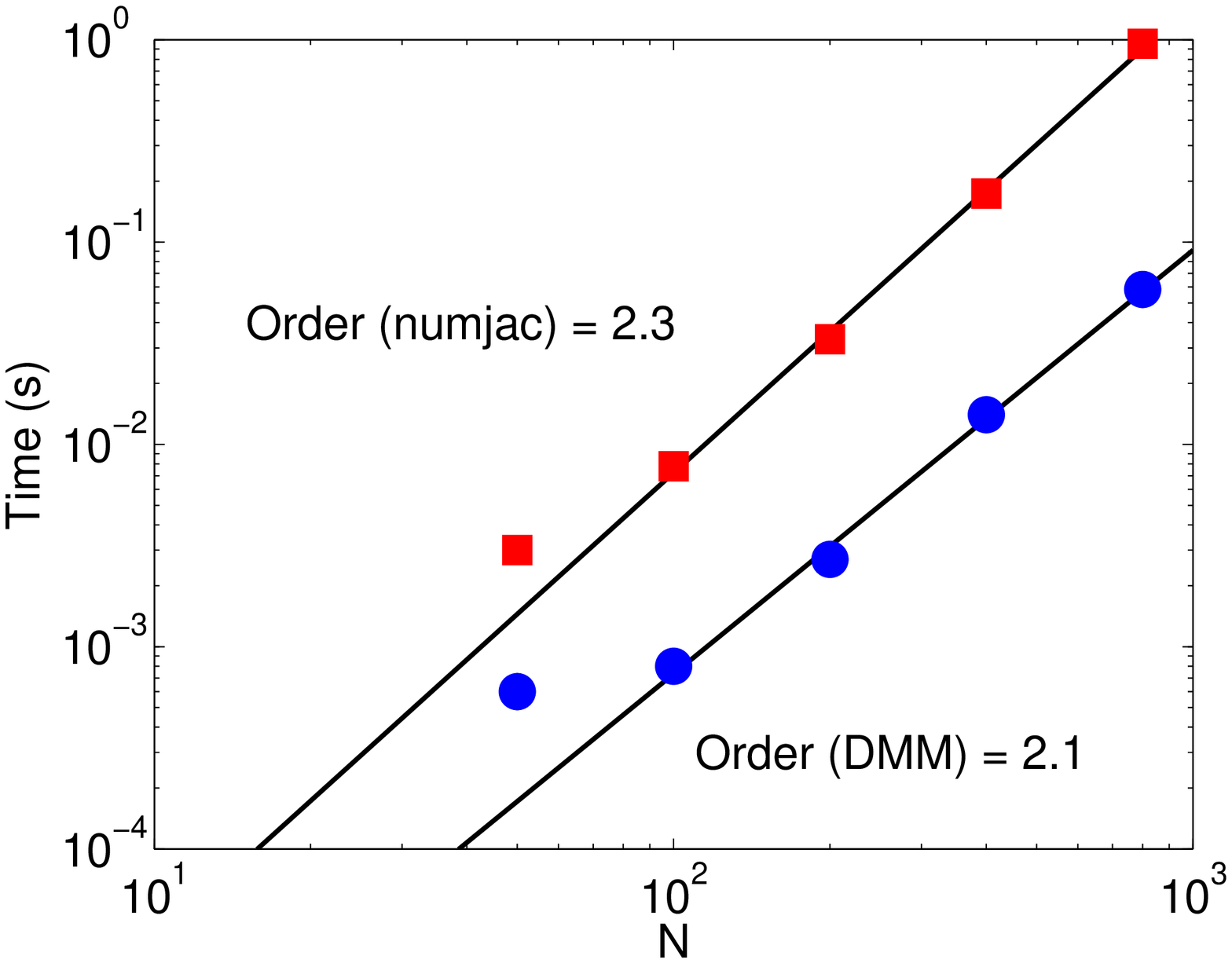}}
\ \ \ \
\scalebox{0.35}{\includegraphics{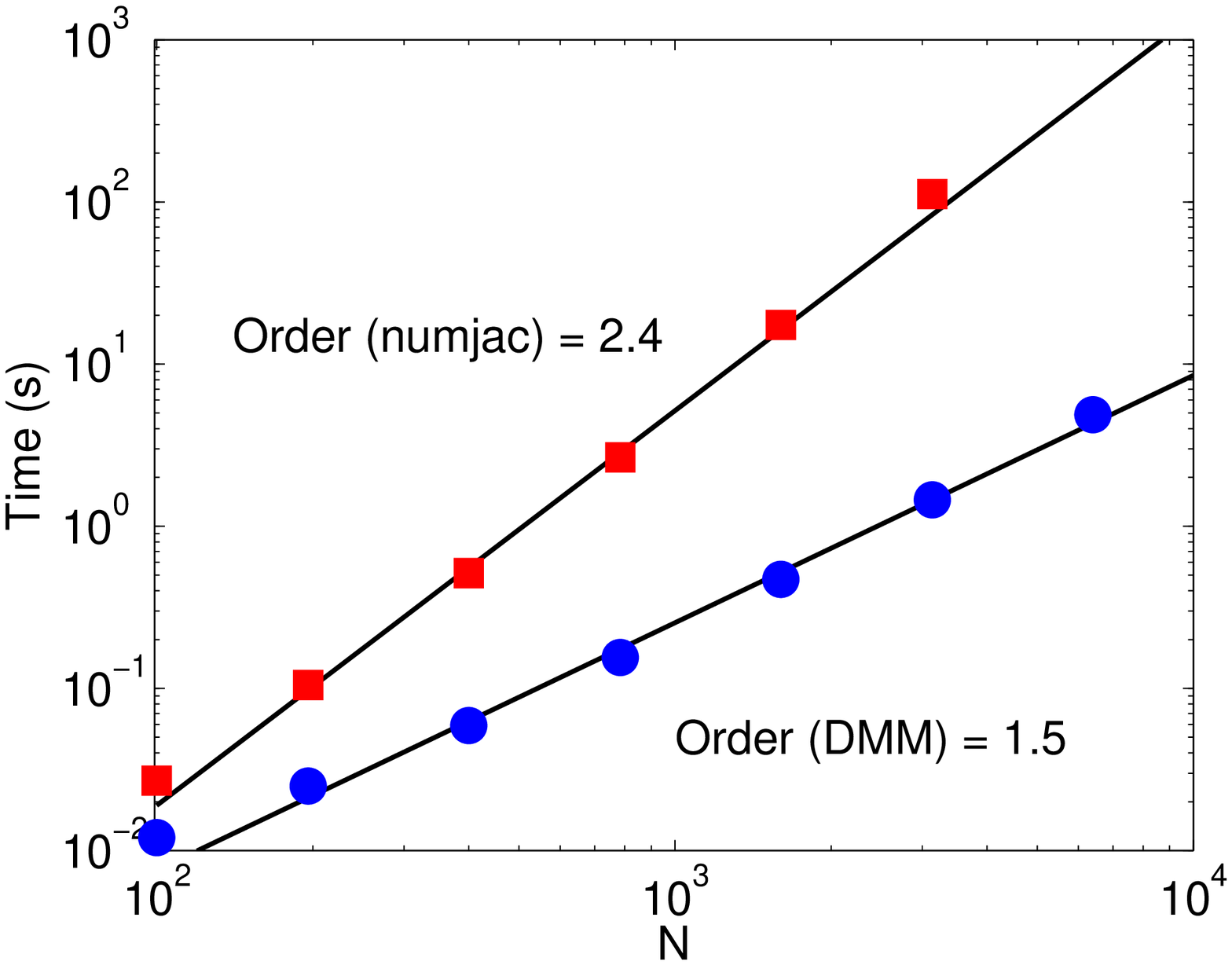}}
\caption[Computational Performance of Sparse vs. Dense Representations]{
\label{figure:jacobian_construction_perf}
Comparison of the computational performance of direct matrix method (circles) 
and numerical Jacobian computed using MATLAB \texttt{numjac()} function 
(squares) for the 1D electrochemical thin-film problem 
(Section~\ref{sec:echem_thin_films}) and the 2D metal colloid sphere problem 
(Section~\ref{sec:dlc_metal_colloid}).  For both comparisons, the
\texttt{numjac()} and direct matrix method codes were optimized by vectorizing
the residual functions and avoiding matrix multiplications whenever
possible.  The data for these graphs were generated on a 2.4 GHz MacBook Pro.
}
\ec
\end{figure}

Figure~\ref{figure:jacobian_construction_perf} compares
the performance of the direct matrix method against the MATLAB 
\texttt{numjac()} function for the two example problems discussed in 
Section~\ref{sec:applications}.  As we can see, the direct matrix method
is at least an order of magnitude faster for both the 1D and 2D problems.
For the 2D problem, the direct matrix method also shows superior scaling 
with the grid size.  To ensure a fair comparison, we vectorized the 
residual calculation to minimize the number of function calls required
by \texttt{numjac()} and avoided the use of matrix 
multiplications\footnote{For example, we express matrix-vector products of 
the form $\diag \left( \hat{u} \right) * \hat{v}$ as component-wise 
multiplication of two grid functions $\hat{u} \dotstar \hat{v}$.},
whenever possible, which benefited both methods.  Matrix-matrix 
multiplications are especially detrimental for the direct matrix method 
because they can worsen the scaling of the Jacobian construction time with 
grid size to the point where the numerical Jacobian is faster to compute.
For instance, in the left graph in 
Figure~\ref{figure:jacobian_construction_perf}, a Jacobian computed using
the direct matrix method with explicit matrix-matrix multiplications
take $O(N^3)$ time, which negates the performance benefits of the 
method compared with a \texttt{numjac()} implementation before $N$ even 
reaches 1000.  

In addition to avoiding matrix-matrix multiplication, it is important to use 
sparse matrices when possible.  For problems in more than one space dimension, 
sparse matrices are produced when Kronecker products with identity matrices 
are used to construct differentiation matrices even if the 1D differentiation
matrices are dense.  Not only does the memory required for dense matrix 
representations easily exhaust the memory on workstations and laptops,  
dense matrix representations also leads to poor computational performance 
when applying and multiplying the matrices.  In general, sparse matrix
operations have better scaling properties as the grid size grows.

\subsection{Relationship to the Newton-Kantorovich Method}
The direct matrix method for computing the Jacobian of discretized
integro-differential equations is closely related to the calculation of
the Fr\'echet derivative\footnote{Recall that the Fr\'echet derivative for 
nonlinear functionals is the generalization of the Jacobian for nonlinear
functions over finite-dimensional 
spaces~\cite{boyd_spectral_book,michel_book,taylor_1974}.  
For intuition, Ortega and Rheinboldt provide a nice discussion of Fr\'echet 
derivatives in the context of finite-dimensional spaces~\cite{ortega_book}.}
used in the Newton-Kantorovich method~\cite{boyd_spectral_book,kelley1987} 
(also known as quasilinearization~\cite{deuflhard_book}). 
The basic idea behind solving nonlinear integro-differential equations using
the Newton-Kantorovich method is to carry out Newton's method in 
\emph{function space}.  For each Newton iteration, we compute the Fr\'echet 
derivative of the integro-differential equation in function space and 
numerically solve the resulting \emph{linear} integro-differential equation 
for the correction to the current iterate of the solution.  Essentially, 
the Newton-Kantorovich method reverses the order of (1)~discretization of 
the continuous problem and (2)~Newton iteration.  Because the equations to be 
solved during each Newton iteration is linear, there is no need to compute a 
Jacobian of the discretized equations.

An important feature of the Newton-Kantorovich method is that the numerical 
discretization used to solve the linearized equation during each Newton 
iteration can, in principle, be completely independent of the discretization 
used to compute the residual of the nonlinear integro-differential equation.  
This freedom can affect the convergence rate of the method depending on the 
degree to which the discretized form of the linearized problem approximates 
the Jacobian of the discretized residual equation.

Because the direct matrix method begins with a discrete equation possessing
the same mathematical structure as the continuous residual equation, it 
produces a Jacobian that is a discrete analogue of the Fr\'echet derivative
for the continuous integro-differential equation.  Unlike the 
Newton-Kantorovich method, however, the direct matrix method produces the
unique Jacobian associated with the particular choice of discretization for 
the residual of the nonlinear integro-differential equation.  The 
freedom to independently choose the numerical discretizations for the 
residual equation and the Fr\'echet derivative is not present in the direct 
matrix method.  As a result, given a numerical discretization for the 
residual equation, the direct matrix method can be viewed as a way to 
generate the optimal discretization for the linearized equation that arises 
during each Newton iteration of the Newton-Kantorovich method.

\section{Applications \label{sec:applications}}
Analysis of electrochemical systems is a classical subject that has recently 
seen a renewal of interest.  Modern electrochemical systems of interest
include ion channels in biological 
membranes~\cite{barcilon1992,barcilon1997,park1997},
microfluidic devices based on electro-osmotic 
flows~\cite{bazant2004,squires2004}, 
and thin-film battery technologies~\cite{neudecker2000,takami2002,wang1996}.
A common feature of many of these applications is that the electrochemical
system is operated under extreme conditions, such as large applied fields
or very small physical size ~\cite{bazant2005,chu2005}.
In these regimes, numerical solutions of the nonlinear governing equations 
are useful for gaining insight into the rich behavior of these systems.
As we shall see, the direct matrix method makes it easy to compute the 
analytical Jacobian required to solve these nonlinear equations using Newton's
method.

\subsection{Electrochemical Thin-Films \label{sec:echem_thin_films}}
Analysis of 1D electrochemical systems leads to an example of a nonlinear 
integro-differential equation.  For steady-state electrochemical thin-films 
made up of a dilute solution of symmetric binary electrolyte with faradaic 
reactions at the surfaces of the thin-film~\cite{bazant2005,chu2005}, the 
electric field, $E$, satisfies the second-order differential 
equation\footnote{This equation is mathematically
equivalent to the Poisson-Nernst-Planck equation formulation of 
electrochemical transport~\cite{bazant2005}.  To simplify the discussion,
equation (\ref{eq:echem_thin_film_eqn}) is a slightly modified form of the 
master equation in~\cite{bazant2005,chu2005} derived by making the
substitutions $x \rightarrow (x+1)/2$ and $E \rightarrow 2E$.}
\beq
\epsilon^2 \left( \frac{d^2 E}{dx^2} - \frac{1}{2} E^3 \right)
- \frac{1}{4} \left(c_0 + j (x+1) \right) E - \frac{j}{4} = 0
\label{eq:echem_thin_film_eqn}
\eeq
on the domain $(-1, 1)$ subject to boundary conditions that
represent the kinetics of electrode reactions
\bea
  - k_c (c(1) + \rho(1)) + j_r - j &=& 0 \\
    k_c (c(-1) + \rho(-1)) - j_r - j &=& 0,
\label{eq:echem_thin_film_bcs}
\eea
where $c$ is the average ion concentration, $\rho$ is the charge density,
$j$ is the current density flowing through the thin-film, $\epsilon$
is a parameter related to the dielectric constant, $k_c$ and $j_r$ are 
reaction rate constants, and $c_0$ is the following expression
\beq
  c_0 = (1-j) 
      + \epsilon^2 \left[2 E(1) - 2 E(-1) - \int_{-1}^1 E^2 dx \right].
\eeq
The average ion concentration and charge density are related to the electric 
field via the equations
\beq
  c(x) = c_0 + j (x+1) + 2 \epsilon^2 E^2 \ , \ 
  \rho(x) = 4 \epsilon^2 \frac{dE}{dx}. 
\eeq

We can solve this set of equations via Newton's method using a systematic 
application of the direct matrix method.  To discretize the equations, we
use a pseudospectral method based on the Chebyshev grid on the interval 
$[-1,1]$.  The differentiation matrix, $D$, for this computational grid is 
just the standard differentiation matrix for the Chebyshev 
grid~\cite{boyd_spectral_book,fornberg_spectral_book,trefethen_spectral_book}.  
For numerical integration, we use the Clenshaw-Curtis 
quadrature weights~\cite{trefethen_spectral_book}, which we denote by the 
row vector $w$.  The quadrature weights are used to construct a quadrature
matrix, $Q$, which is the analogue of the differentiation matrix: 
$Q = [w^T, w^T, \ldots, w^T]^T$.  When a grid function $f$ is multiplied by 
$Q$, the result is a vector where all entries are equal to the numerical
approximation of the integral of $f$.  

With these discrete operators, we can put (\ref{eq:echem_thin_film_eqn})
into matrix-vector form:
\beq
\epsilon^2 \left( D^2 * \hat{E} - \frac{1}{2} \hat{E} \dotwedge 3 \right) 
- \frac{1}{4}\left(\hat{C}_0 + j (x+1) \right) .* \hat{E} - \frac{j}{4} = 0
  \label{eq:discrete_echem_thin_film_eqn}
\eeq
with 
\beq
\hat{C}_0 = (1-j) + \epsilon^2 \left( 2 \hat{E}_1 - 2 \hat{E}_N
                  - Q * (\hat{E} \dotwedge 2) \right),
\eeq
where we have chosen to order the indices so that $x_1 = 1$ and $x_N = -1$
(this follows the convention used in~\cite{trefethen_spectral_book} and in
the code in Appendix \ref{appendix:echem_thin_film_code}).  
The boundary conditions are imposed by replacing the discrete equations
corresponding to $x_1$ and $x_N$ with
\bea
  -k_c \left( \hat{c}_0 + 2 j + \epsilon^2 
      \left( 2 \hat{E}_1^2 + 4 D_1 * \hat{E} \right) 
      \right) + j_r - j = 0 \\
  k_c \left( \hat{c}_0 + \epsilon^2 
      \left( 2 \hat{E}_N^2 + 4 D_N * \hat{E} \right) 
      \right) - j_r - j = 0 ,
  \label{eq:discrete_echem_thin_film_bc}
\eea
where $D_1$ and $D_N$ are the rows of the differentiation matrix corresponding
to $x_1$ and $x_N$, respectively, and $\hat{c}_0$ is a single component of 
$\hat{C}_0$.
\begin{figure}[tb] 
\bc
\scalebox{0.35}{\includegraphics{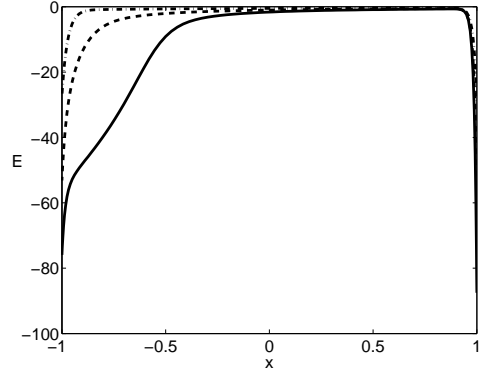}} 
\caption{
\label{fig:echem_thin_film_solns}
Solution of electrochemical thin-film equations (\ref{eq:echem_thin_film_eqn}) 
-- (\ref{eq:echem_thin_film_bcs}) computed using $100$ grid points with 
$\epsilon = 0.01$, $k_c = 10$, and $j_r = 10$ for $j = 1.5$ (solid),
$j = 1.0$ (dash), and $j = 0.5$ (dot-dash).
}
\ec
\end{figure}

The Jacobians for the left-hand side of these discrete equations are now 
easily computed.  Applying the differentiation rules from 
Section~\ref{sec:direct_matrix_method}, the Jacobian for the interior
grid points is
\beq
  J_{int} = 
  \epsilon^2 \left( D^2 - \frac{3}{2} 
    \diag\left( \hat{E} \dotwedge 2 \right)  \right)
  - \frac{1}{4} \diag\left(\hat{C}_0 + j (x+1)\right) 
  - \frac{1}{4} \diag\left(\hat{E}\right) 
    * \frac{\partial \hat{C}_0}{\partial \hat{E}}
  \label{eq:echem_thin_film_J_interior}
\eeq
with 
\beq
  \frac{\partial \hat{C}_0}{\partial \hat{E}} = 
  \epsilon^2 \left( 
  \left[ 
  \begin{array}{ccccc}
  2 & 0 & \cdots & 0 & -2 \\
  \vdots & \vdots & \vdots & \vdots & \vdots \\
  2 & 0 & \cdots & 0 & -2 
  \end{array}
  \right] 
  - 2 Q*\diag\left(\hat{E}\right) \right)
\eeq
The Jacobian for the discretized boundary conditions are similarly calculated:
\bea
  J_1 &=&  -k_c \left( \frac{\partial \hat{c}_0}{\partial \hat{E}} + 
           4 \epsilon^2 [E_1 \ 0 \ \ldots \ 0 ] + 4 \epsilon^2 D_1 \right)  \\
  J_N &=&  k_c \left( \frac{\partial \hat{c}_0}{\partial \hat{E}} + 
           4 \epsilon^2 [0 \ \ldots \ 0 \ E_N] + 4 \epsilon^2 D_N \right) 
  \label{eq:echem_thin_film_J_bcs},
\eea
where
\beq
  \frac{\partial \hat{c}_0}{\partial \hat{E}} =
  \epsilon^2 \left( [2 \ 0 \ \ldots \ 0 \ -2]
                  - 2 w*\diag\left(\hat{E}\right) \right) 
\eeq

From the perspective of computational performance, the above formulation 
of the Jacobian is suboptimal because it includes a matrix-matrix multiply
in (\ref{eq:echem_thin_film_J_interior}) that can be avoided.  
To reduce the time required to compute the Jacobian, the key observation 
is that each row of $\frac{\partial \hat{C}_0}{\partial \hat{E}}$ is 
equal $\frac{\partial \hat{c}_0}{\partial \hat{E}}$.  Therefore, 
$\diag\left(\hat{E}\right) * \frac{\partial \hat{C}_0}{\partial \hat{E}}$ 
is more efficiently computed as the Kronecker product of 
$\hat{E}$ and $\frac{\partial \hat{c}_0}{\partial \hat{E}}$.  The 
evaluation of the residual can also be improved by recognizing that all of
the elements of $\hat{C}_0$ are equal to $\hat{c}_0$, but this optimization
has a far smaller impact than the reformulation of the Jacobian.

\begin{figure}[tb] 
\bc
\scalebox{0.35}{\includegraphics{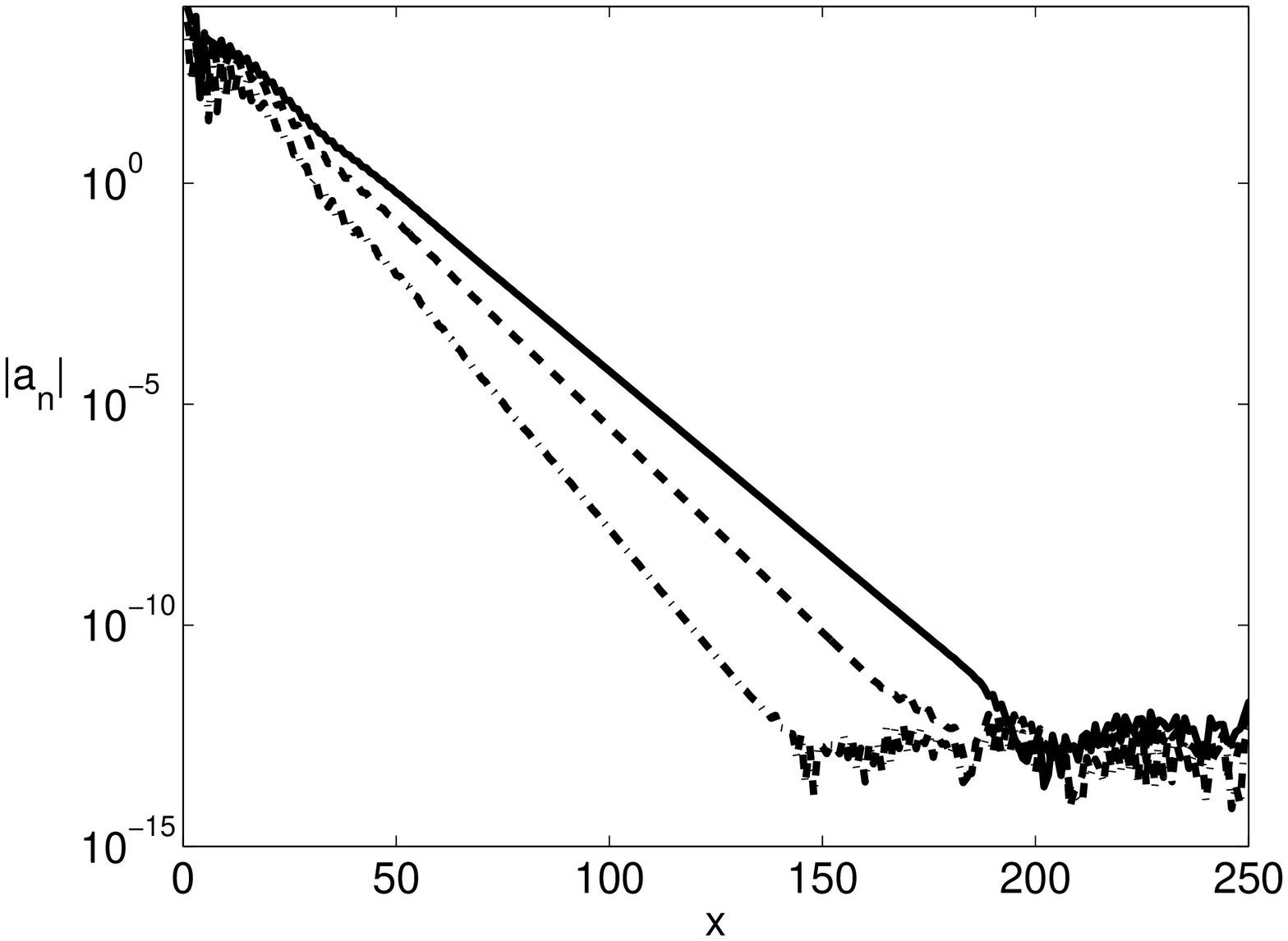}} 
\ \ \ \
\scalebox{0.35}{\includegraphics{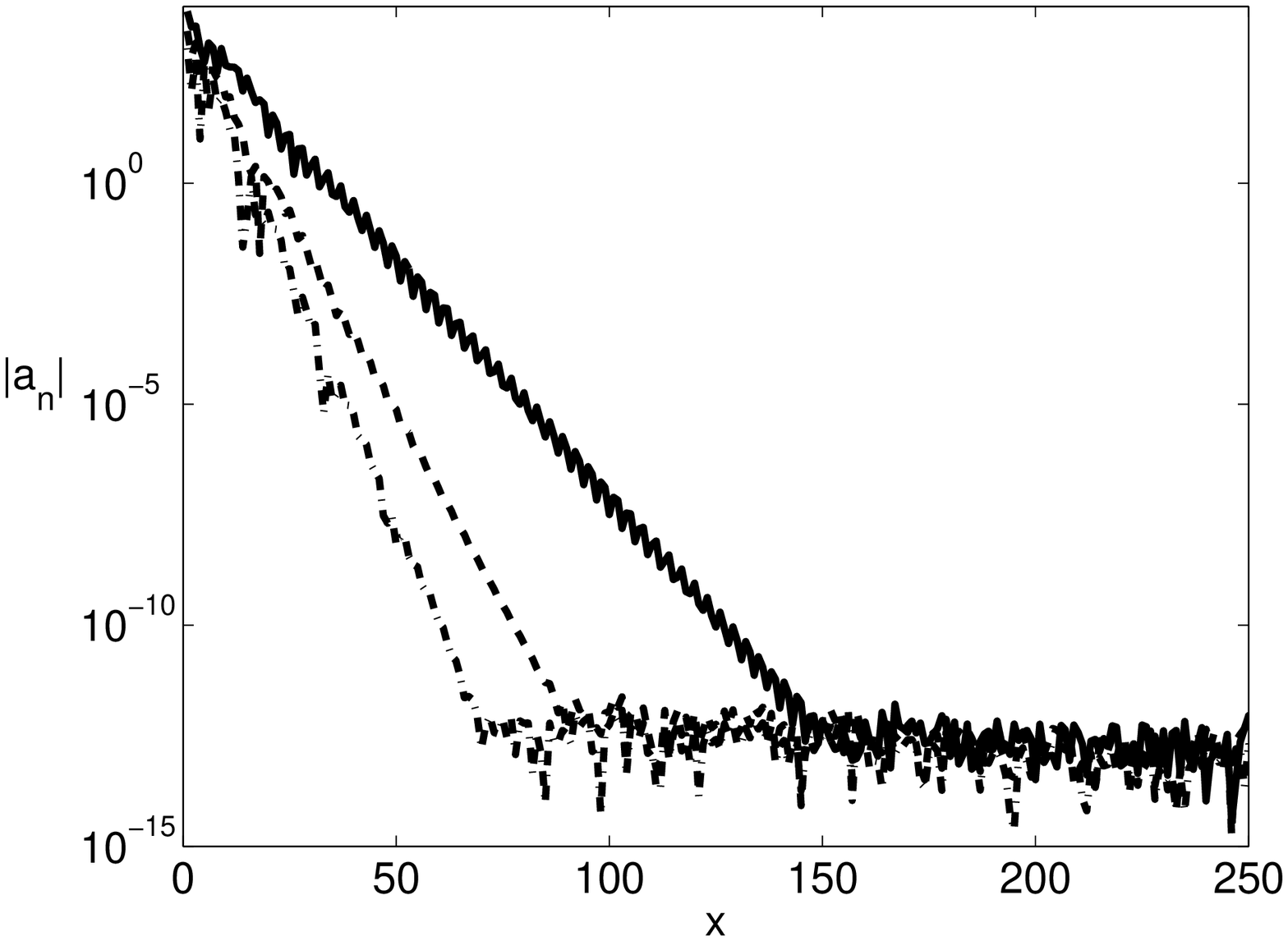}} 
\caption{
\label{fig:echem_thin_film_convergence}
Plots of the absolute value of the spectral coefficients, $a_n$, for the 
numerical solution of the electrochemical thin-film equations with 
$j=1.5$ (solid), $j=1.0$ (dash), and $j=0.5$ (dot-dash) as a function
of the basis function degree when (\ref{eq:echem_thin_film_eqn}) -- 
(\ref{eq:echem_thin_film_bcs}) are solved directly (left) and using the
variable transformation (\ref{eq:echem_var_transform}) to place more grid 
points in the boundary layers (right).  Although the convergence rate as a 
function of number of grid points is geometric in all cases, the variable 
transformation significantly reduces the number of grid points required obtain 
the best solution possible given the precision of the computation.
Note that the spectral coefficients were calculated from a numerical solution 
generated using $250$ grid points to make the roundoff plateau more apparent.
}
\ec
\end{figure}

Now that we have explicitly computed all of the components for Newton's 
method, it is straightforward to write a program to solve the
nonlinear integro-differential equations for electrochemical thin-films.
The MATLAB code for solving is relatively short and runs quickly (see 
Appendix \ref{appendix:echem_thin_film_code}).  
One special issue that arises for this problem is that continuation 
methods~\cite{boyd_spectral_book} are required to obtain good initial 
iterates for the Newton iteration at high current densities.  
Figure~\ref{fig:echem_thin_film_solns} shows the numerical solution of
(\ref{eq:echem_thin_film_eqn}) -- (\ref{eq:echem_thin_film_bcs}) computed
using $100$ grid points with $\epsilon = 0.01$, $k_c = 10$, and $j_r = 10$ 
for various values of $j$.  As expected, we observe geometric convergence
with respect to the number of grid points (see 
Figure~\ref{fig:echem_thin_film_convergence}).  Notice that at higher
current densities, we see slower convergence rates due to the presence of
greater structure in the solution.  

While quite satisfactory, the convergence rate for the numerical 
discretization (\ref{eq:discrete_echem_thin_film_eqn}) -- 
(\ref{eq:discrete_echem_thin_film_bc}) as a function of the number of 
grid points is limited by the need to resolve the boundary layers.
By using a mapping of the computational domain that allows us to 
place a few grid points within the boundary layers, we can obtain 
a faster convergence rate.  For example, by using the variable
transformation: 
\bea
x = \frac{1}{\beta}\tanh\left(\mathrm{atanh}(\beta) y\right) \ , \ 
E(x) = \frac{\beta}{\mathrm{atanh}\beta} 
         \cosh^2 \left(\mathrm{atanh}(\beta) y \right) E(y),
\label{eq:echem_var_transform}
\eea
where $\beta$ is an adjustable parameter less than $1$, we can 
significantly reduce the number of grid points required obtain 
a solution accurate to machine precision (see 
Figure~\ref{fig:echem_thin_film_convergence}).  It is interesting to 
note that the optimal value for $\beta$ depends on the current density
$j$.  For $j=0.5$ and $j=1.0$, a $\beta$ value of $0.9$ yields near
optimal results.  For $j=1.5$, however, the fastest convergence is obtained
near $\beta=0.75$.
As is typical, the transformed electrochemical thin-film equations are bit 
more complicated to deal with than the original equations.  However, the 
direct matrix method makes it straightforward to discretize the transformed 
equations and compute the exact Jacobian for the resulting nonlinear algebraic 
equations (see Appendix \ref{appendix:echem_thin_film_mapped_code}).

\subsection{Double Layer Charging of Metal Colloid Sphere at 
            High Applied Electric Fields \label{sec:dlc_metal_colloid}}
Analysis of double layer charging for colloid systems subject to applied 
electric fields gives rise to nonlinear differential equations in multiple
space dimensions with complicated boundary conditions.  In the electroneutral 
limit~\cite{bazant2005,newman_book,rubinstein_book}, the steady-state 
governing equations for systems composed of symmetric binary electrolyte 
are~\cite{chu2006}
\bea 
  \lapl c &=& 0 
  \label{eq:c_eqn_steady}
  \\
  \div \left( c \grad \phi \right)  &=& 0, 
  \label{eq:phi_eqn_steady}
\eea
where $c$ is the average ion concentration and $\phi$ is the electric 
potential.  For metal colloid surfaces, the appropriate 
boundary conditions are~\cite{chu2006,chu2007}
\bea
  0 &=& \epsilon \divs \left( q \grads \ln c + w \grads \phi \right) 
      - c \frac{\partial \phi}{\partial n} \\
  0 &=& \epsilon \divs \left( w \grads \ln c + q \grads \phi \right) 
      - \frac{\partial c}{\partial n} \\
  q &=& -2 \sqrt{c} \sinh(\zeta/2) 
  \label{eq:dlc_q_def_GCS} \\
  w &=& 4 \sqrt{c} \sinh^2(\zeta/4) 
  \label{eq:dlc_w_def_GCS} \\
  v - \phi &=& \zeta + 2 \delta \sqrt{c} \sinh (\zeta/2) 
  \label{eq:dlc_stern_bc_GCS}
\eea
where $q$ and $w$ are the excess charge and ion concentration in the 
boundary layer, $\zeta$ is the electric potential drop across the boundary
layer, $v$ is the potential of the metal colloid, and $\delta$ is a parameter
related to the capacitance of the boundary layer.

As a model problem, we solve these equations for a metal colloid sphere 
subjected to a uniform applied electric field of strength $E$ in the 
$z$-direction.
To avoid infinite values of the electric potential, the numerical model is 
formulated in terms of $\psi \equiv \phi + Ez$, the deviation of 
the electric potential from that of the uniform applied field.  The
spherical geometry of the problem also allows us to demonstrate the use of 
the direct matrix method on a non-Cartesian (though still logically 
rectangular) grid.

While this problem may seem daunting, it is straightforward to obtain a
solution numerically by using Newton's method with an analytical Jacobian
computed using the direct matrix method.  Taking advantage of azimuthal 
symmetry, we discretize the equations in spherical coordinates on a 2D 
pseudospectral grid that is the tensor product of grids in the radial and 
polar angle directions.  We use a shifted semi-infinite rational Chebyshev 
grid~\cite{boyd_spectral_book} in the radial direction and a uniformly 
spaced grid for the polar angle direction.  The required differentiation 
matrices are constructed using Kronecker products, and the boundary conditions 
are handled using restriction and prolongation matrices as discussed in 
Section~\ref{sec:multiple_space_dims}.  

To facilitate the formulation of the matrix-vector representation of the 
equations, let us fix our notation.  Let $D_r$ and $D_\theta$ be the radial 
and angular contributions to the discrete divergence operator, $G_r$ and 
$G_\theta$ be the radial and angular components of the discrete gradient
operator, and $L$ be the discrete Laplacian operator.  Also, let $n$ and $s$ 
subscripts denote normal and tangential derivative operators at the surface
of the sphere.  

For the purpose of discussion (and implementation), it is convenient to 
decompose the discrete differential operators into pieces that correspond to
contributions from finite and infinite grid points.  For example, $L$ can be 
decomposed into $L^f$ and $L^\infty$ which respectively account for the 
contributions to the Laplacian operator from finite and infinite grid points; 
that is, $L*\ch = L^f * \ch_f + L^\infty * \ch_\infty$, where $\ch_f$ and 
$\ch_\infty$ are the concentration values at finite and infinite grid points 
respectively.  Similarly, to impose the boundary conditions, we use derivative 
operators that act only on surface values.  Surface operators and surface 
field values will be denoted with superscripts $s$ and subscripts $s$, 
respectively.  
Finally, to refer to values at interior grid points (\ie finite grid points 
that are \emph{not} on the surface of the sphere), we use the subscript $i$. 

In this notation, the discretized form of the bulk equations
(\ref{eq:c_eqn_steady}) and (\ref{eq:phi_eqn_steady}) are given by
\bea
  0 = F_1 &\equiv& L^f * \ch_f + L^\infty * c_\infty  
  \label{eq:num_mod_diffusion_steady} \\
  0 = F_2 &\equiv& 
          D_r^f * \left[ \ch_f \dotstar
               \left( G_r^f * \psih_f - E \cos \theta \right) \right]
        - D_r^\infty * \left( c_\infty \dotstar E \cos \theta \right)
  \nonumber \\
        &+& 
          D_\theta^f * \left[ \ch_f \dotstar 
             \left( G_\theta^f * \psih_f + E \sin \theta \right) \right]
        + D_\theta^\infty * \left( c_\infty \dotstar E \sin \theta \right).
  \label{eq:num_mod_current_steady} 
\eea
In these equations, the unknowns are the values of the $c$ and $\psi$ at 
finite grid points; values at infinity are specified by the boundary
conditions and so are known quantities (which is why $c_\infty$ does not
have a hat accent and $\psi_\infty = 0$ does not show up at all). 
In discretized form, the boundary conditions on the surface of the sphere
are 
\bea
  0 = H_1 &\equiv& 
   \epsilon D_s * \left[ 
          \qh \dotstar \left( G^s * \ln \ch_s \right) 
        + \wh \dotstar \left( G^s * \psih_s - G^s * E \cos \theta \right) 
       \right]
   \nonumber \\ 
   &-& c_s \dotstar \left( G_n^f * \psih_f + E \cos \theta \right) 
  \label{eq:num_mod_bc_q_steady} \\
  0 = H_2 &\equiv& 
   \epsilon D_s * \left[ 
         \wh \dotstar \left( G^s * \ln \ch_s \right) 
       + \qh \dotstar \left( G^s * \psih_s - G^s * E \cos \theta \right) 
       \right]
   \nonumber \\ 
   &-& \left( G_n^f * \ch_f + G_n^\infty * c_\infty \right).
  \label{eq:num_mod_bc_w_steady}
\eea
Closure for these equations is given by using 
(\ref{eq:dlc_q_def_GCS}) -- (\ref{eq:dlc_w_def_GCS}) to relate
$\qh$ and $\wh$ to the zeta-potential and using (\ref{eq:dlc_stern_bc_GCS})
to compute the zeta-potential from $\phi$ and $\ch_s$.

The direct matrix method makes it straightforward to derive the analytical
Jacobian for the system of equations
(\ref{eq:dlc_q_def_GCS}) -- (\ref{eq:num_mod_bc_w_steady}).
The derivatives of $F_1$ and $F_2$ with respect to the unknowns $\ch_f
$ and $\psih_f$ are easily calculated:
\bea
  \DDx{F_1}{\ch_f} &=& L^f \\
  \DDx{F_1}{\psih_f} &=& 0 \\
  \DDx{F_2}{\ch_f} &=& 
    D_r^f * \diag \left( G_r^f * \psih_f - E \cos \theta \right) 
   +D_\theta^f * \diag \left( G_r^f * \psih_f + E \sin \theta \right) 
    \\
  \DDx{F_2}{\psih_f} &=& 
    D_r^f * \diag \left( \ch_f \right) * G_r^f
   +D_\theta^f * \diag \left( \ch_f \right) * G_\theta^f
\eea
The derivatives for the discretized boundary conditions are more
complicated because $\qh$, $\wh$, and $\ch_s$ implicitly depend on 
the unknown variables and because surface grid points must be treated 
differently than interior grid points.  However, a systematic application of 
the differentiation rules in Section~\ref{sec:differentiation_rules} yields
the analytical Jacobian directly in matrix form:
\bea
  \DDx{H_1}{\ch_s} &=& 
     \frac{\epsilon}{2} D_s * 
     \diag \left(~q \dotslash c_s 
                \dotstar \left ( G^s * \ln c_s \right) ~\right)
  \nonumber \\
    &-& \epsilon D_s * \diag \left(
        \sqrt{c_s} \dotstar \cosh(\zeta/2) 
                   \dotstar \frac{\partial \zeta}{\partial c_s} 
                   \dotstar \left ( G^s * \ln c_s \right) ~\right)
  \nonumber \\
    &+& \epsilon D_s * \diag (q) * G^s *
                 \diag \left(1 \dotslash c_s\right)
  \nonumber \\
    &+& \frac{\epsilon}{2} D_s * \diag \left(~
      w \dotslash c_s
      \dotstar \left ( G^s * \psi_s - G^s * E \cos \theta \right)~\right)
  \nonumber \\
    &+& \epsilon D_s *\diag \left(
        \sqrt{c_s} \dotstar \sinh(\zeta/2)  
        \dotstar \frac{\partial \zeta}{\partial c_s}
        \dotstar \left( G^s * \psi_s - G^s * E \cos \theta \right) 
        \right)
  \nonumber \\
  &-& \diag \left( G_n^f * \psi_f + E \cos \theta \right) 
  \label{eq:dlc_jacobian_H1}
  \\
  \DDx{H_1}{c_i} &=& 0 \\
  \DDx{H_1}{\psi_s} &=& 
    -\epsilon D_s * 
      \diag \left(~\sqrt{c_s} \dotstar \cosh (\zeta/2) 
      \dotstar \frac{\partial \zeta}{\partial \psi_s} 
      \dotstar (G^s * \ln c_s)~\right)
    + \epsilon D_s * \diag(w) * G^s
  \nonumber \\ 
    &+& \epsilon D_s * \diag \left(
        \sqrt{c_s} \dotstar \sinh \left( \zeta/2 \right) 
        \dotstar \frac{\partial \zeta}{\partial \psi_s} 
        \dotstar \left( G^s * \psi_s - G^s * E \cos \theta \right)
        \right)
  \nonumber \\ 
  &-& \diag \left( c_s \right) * G_n^s
  \\
  \DDx{H_1}{\psi_i} &=& -\diag \left( c_s \right) G_n^i
\eea
\bea
  \DDx{H_2}{\ch_s} &=& 
     \frac{\epsilon}{2} D_s * 
     \diag \left(~w \dotslash c_s 
                \dotstar \left ( G^s * \ln c_s \right) ~\right)
    \nonumber \\
    &+& \epsilon D_s * \diag \left(
      \sqrt{c_s} \dotstar \sinh(\zeta/2) 
                 \dotstar \frac{\partial \zeta}{\partial c_s} 
                 \dotstar \left ( G^s * \ln c_s \right) ~\right)
  \nonumber \\
    &+& \epsilon D_s * \diag (w) * G^s *
                 \diag \left(1 \dotslash c_s\right) 
    \nonumber \\
    &+& \frac{\epsilon}{2} D_s * \diag \left(~
      q \dotslash c_s
      \dotstar \left ( G^s * \psi_s - G^s * E \cos \theta \right)~\right)
  \nonumber \\ &-&
      \epsilon D_s *\diag \left(
       \sqrt{c_s} \dotstar \cosh(\zeta/2)  
       \dotstar \frac{\partial \zeta}{\partial c_s}
       \dotstar \left( G^s * \psi_s - G^s * E \cos \theta \right) 
       \right)
  \nonumber \\
  &-& G_n^s
  \\
  \DDx{H_2}{c_i} &=& -G_n^i \\
  \DDx{H_2}{\psi_s} &=& 
      \epsilon D_s * 
      \diag \left(~\sqrt{c_s} \dotstar \sinh (\zeta/2) 
      \dotstar \frac{\partial \zeta}{\partial \psi_s} 
      \dotstar (G^s * \ln c_s)~\right)
    + \epsilon D_s * \diag(q) * G^s 
  \nonumber \\ 
    &-& \epsilon D_s * \diag \left(
        \sqrt{c_s} \dotstar \cosh \left( \zeta/2 \right) 
        \dotstar \frac{\partial \zeta}{\partial \psi_s} 
        \dotstar \left( G^s \psi_s - G^s * E \cos \theta \right)
        \right)
  \\
  \DDx{H_2}{\psi_i} &=& 0
  \label{eq:dlc_jacobian_H2}
\eea
where
\bea
  \frac{\partial \zeta}{\partial \psi_s} &=& 
    -\frac{1}{1+\delta \sqrt{c_s} \cosh \left( \zeta/2 \right) }
  \\
  \frac{\partial \zeta}{\partial c_s} &=&
    -\frac{\delta \sinh \left( \zeta/2 \right) }
          {\sqrt{c_s} 
           \left [ 1+\delta \sqrt{c_s} \cosh \left( \zeta/2 \right) \right]}.
\eea
The Jacobian for the system of equations is obtained by assembling these
pieces:
\beq
  J = \left[
    \begin{array}{cc}
    \frac{\partial F_1}{\partial \ch} & \frac{\partial F_1}{\partial \psih} \\
    \frac{\partial F_2}{\partial \ch} & \frac{\partial F_2}{\partial \psih} \\
    \frac{\partial H_1}{\partial \ch} & \frac{\partial H_1}{\partial \psih} \\
    \frac{\partial H_2}{\partial \ch} & \frac{\partial H_2}{\partial \psih} 
    \end{array}
    \right],
\eeq
where the Jacobians for $H_1$ and $H_2$ are constructed from 
(\ref{eq:dlc_jacobian_H1}) -- (\ref{eq:dlc_jacobian_H2}) 
using restriction operators.  For instance,
\bea
    \frac{\partial H_1}{\partial \ch} &=& 
      \frac{\partial H_1}{\partial \ch_s} * R_s 
    + \frac{\partial H_1}{\partial \ch_i} * R_i,
\eea
where $R_s$ and $R_i$ are restriction operators for surface and interior
grid points, respectively.  While the formulas may look complicated to 
program, they are actually quite easy to implement in MATLAB (see
Appendix \ref{appendix:metal_colloid_DL_charging_code}). 

\begin{figure}[tb] 
\bc
\scalebox{0.35}{\includegraphics{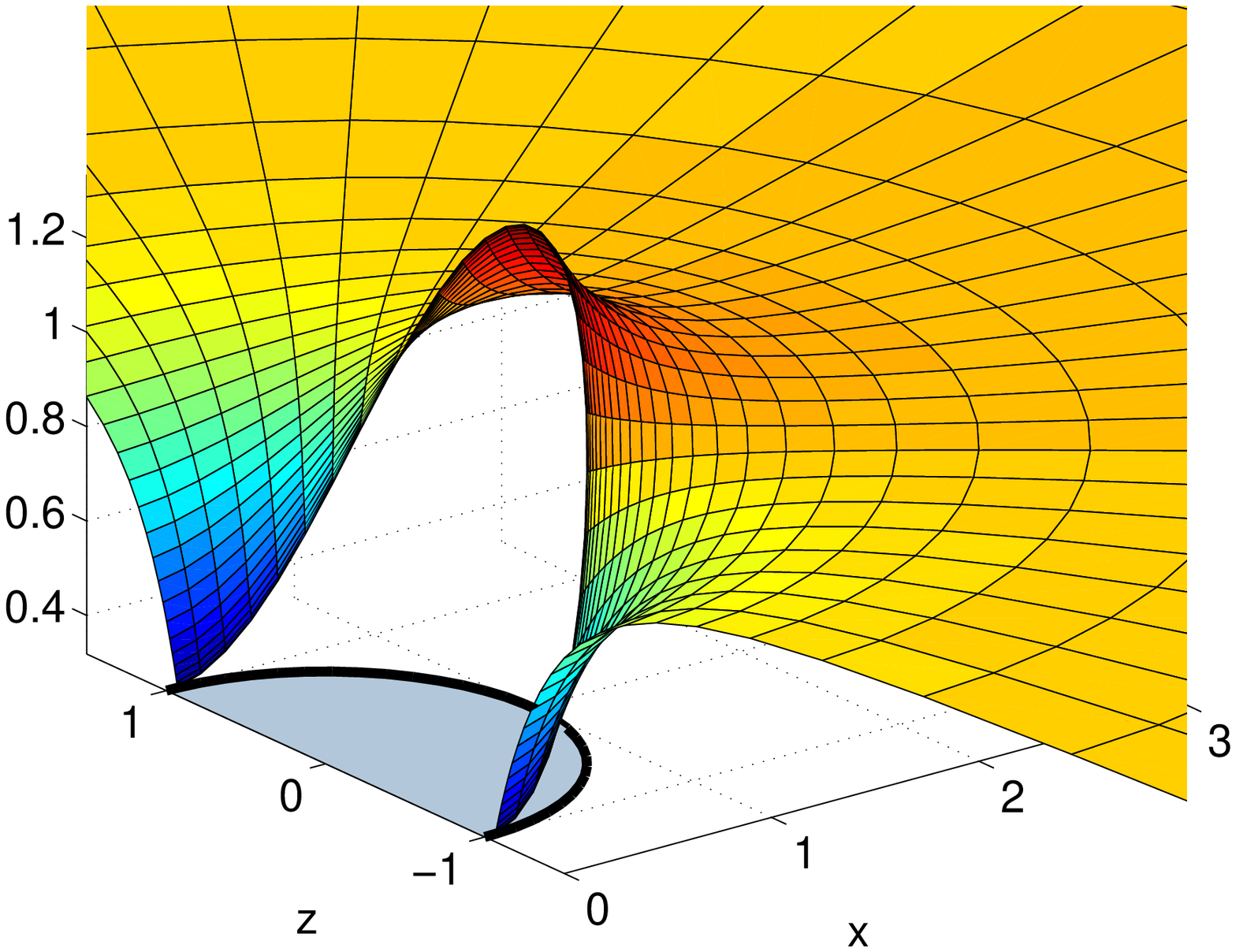}}
\ \ \ \
\scalebox{0.35}{\includegraphics{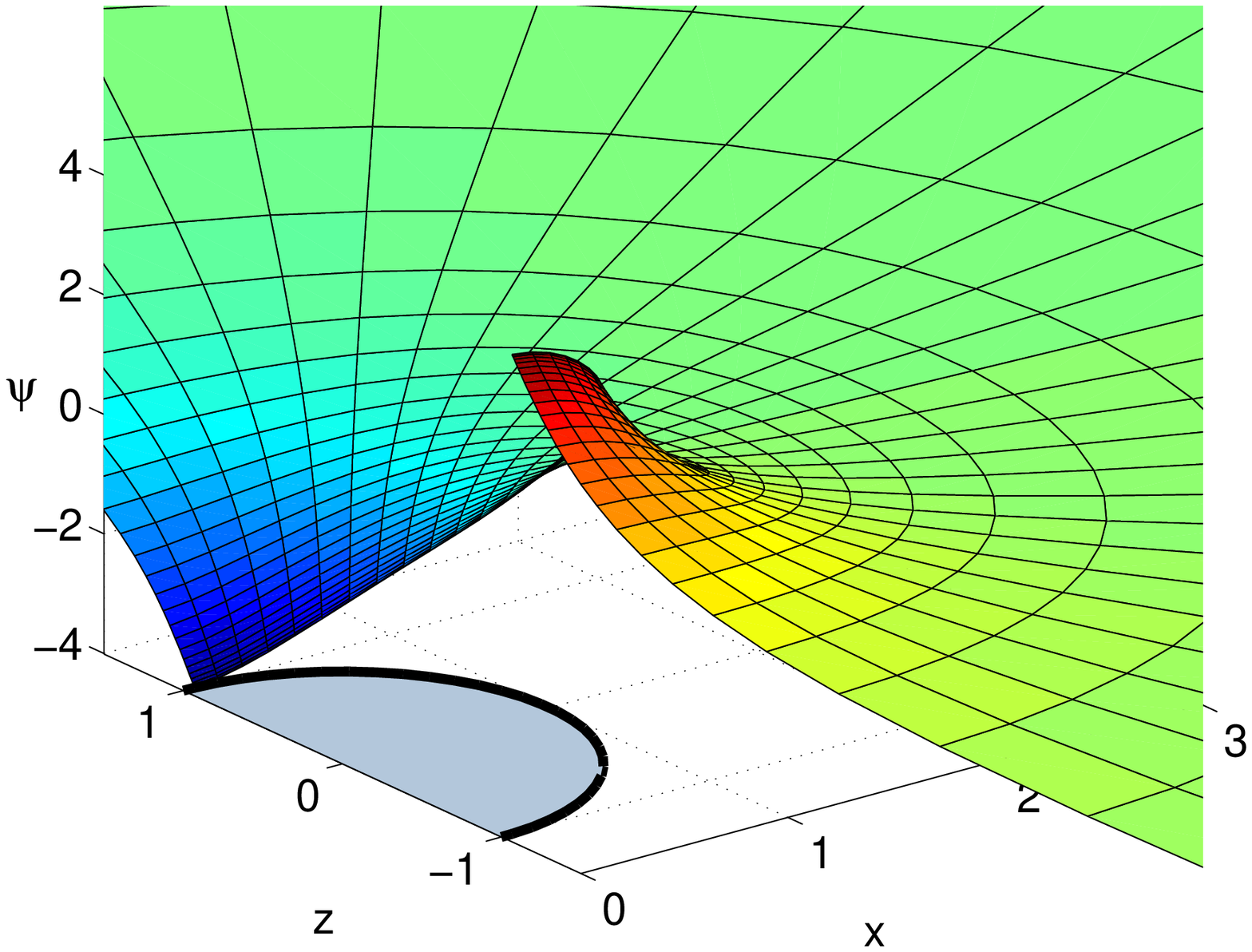}}
\caption{
\label{fig:metal_colloid_DL_charging}
Solution of equations (\ref{eq:c_eqn_steady}) -- (\ref{eq:dlc_stern_bc_GCS})
computed using 30 grid points in both the radial and polar angle directions
with $E = 10$, $v = 0$, $\epsilon = 0.01$, and $\delta = 1$.
}
\ec
\end{figure}

Figure~\ref{fig:metal_colloid_DL_charging} shows numerical solutions obtained 
using the above residual and Jacobian formulas.  As for the electrochemical
thin-film example, continuation is required to obtain good initial iterates
for the Newton iteration at high values of the applied electric fields.  The 
solutions shown are computed for $E = 10$, $v = 0$, $\epsilon = 0.01$,
and $\delta = 1$ using $30$ grid points in both the radial and polar angle 
directions with scale parameter set to $0.5$ for the shifted rational 
Chebyshev grid.  Using pseudospectral grids and the analytical Jacobian, the 
solution is obtained very quickly, requiring only a few Newton iterations for 
each continuation stage (and less than a minute of computation time on a
2.4 GHz MacBook Pro).

\section{Conclusions}
In this article, we have presented a direct matrix method for calculating
analytical Jacobians for discretized, nonlinear integro-differential 
equations.  Because this method is based on simple matrix-based 
differentiation rules, it is less tedious and less error prone than 
other approaches for computing analytical Jacobians.  Furthermore, because 
it yields the Jacobian in matrix form, it is very easy to use languages that 
support vectorized computation to implement numerical methods that require the 
Jacobian.

One interesting possibility that the direct matrix method presents is 
development of high-level automatic differentiation tool for 
discretized nonlinear integro-differential equations.  In contrast to 
traditional automatic differentiation methods~\cite{griewank_book,autodiff_anl}
which operate at the level of individual scalar operations, automatic 
differentiation methods based on the direct matrix method would operate 
on the discrete differential operators associated with the continuous 
differential equation.  Such an automatic differentiation tool could be 
useful for completely eliminating the need for a researcher to compute
the Jacobian of discretized nonlinear integro-differential equations by
hand.

\section*{Acknowledgments}
The author gratefully acknowledges the support of the Department of Energy 
through the Computational Science Graduate Fellowship (CSGF) Program 
provided under grant number DE-FG02-97ER25308, Vitamin D, Inc.,
and the Institute for High-Performance Computing (IHPC) in Singapore.
The author thanks B.~Kim, P.~Fok, and J.~P.~Boyd for many helpful discussions 
and suggestions.

\appendix
\begin{landscape}
\begin{multicols}{2}
\setlength{\columnseprule}{1pt}

\section{\footnotesize MATLAB Code for Electrochemical Thin-Film Example
         \label{appendix:echem_thin_film_code}}
\footnotesize
This code relies on {\texttt cheb.m} and 
{\texttt clencurt.m}~\cite{trefethen_spectral_book}.

\lstset{basicstyle=\tiny}
\begin{lstlisting}
%%%%%%%%%%%%%%%%%%%%%%%%%%%%%%%%%%%%%%%%%%%%%%%%%%%%%%%%%%%%%%%%%%
% parameters
%%%%%%%%%%%%%%%%%%%%%%%%%%%%%%%%%%%%%%%%%%%%%%%%%%%%%%%%%%%%%%%%%%
N = 200;
j = 1.5; epsilon = 0.01; k_c = 10; j_r = 10;    
res_tol = 1e-8; max_iters = 20;

%%%%%%%%%%%%%%%%%%%%%%%%%%%%%%%%%%%%%%%%%%%%%%%%%%%%%%%%%%%%%%%%%%
% compute grid, differentiation matrix, and quadrature weights
%%%%%%%%%%%%%%%%%%%%%%%%%%%%%%%%%%%%%%%%%%%%%%%%%%%%%%%%%%%%%%%%%%
[D,x] = cheb(N-1);         % Chebyshev differentiation matrix 
L     = D*D;               % Laplacian operator
[x,w] = clencurt(N-1);     % Clenshaw-Curtis quadrature weights 

%%%%%%%%%%%%%%%%%%%%%%%%%%%%%%%%%%%%%%%%%%%%%%%%%%%%%%%%%%%%%%%%%%
% set up continuation in j
%%%%%%%%%%%%%%%%%%%%%%%%%%%%%%%%%%%%%%%%%%%%%%%%%%%%%%%%%%%%%%%%%%
j_start = 0.5; dj = 0.1; j_cur = j_start;  

%%%%%%%%%%%%%%%%%%%%%%%%%%%%%%%%%%%%%%%%%%%%%%%%%%%%%%%%%%%%%%%%%%
% generate initial iterate for Newton iteration
%%%%%%%%%%%%%%%%%%%%%%%%%%%%%%%%%%%%%%%%%%%%%%%%%%%%%%%%%%%%%%%%%%
c0 = 1-j_cur; c = c0 + j_cur*(x+1); E = -2*j_cur./(j_cur*(x+1)+c0);

%%%%%%%%%%%%%%%%%%%%%%%%%%%%%%%%%%%%%%%%%%%%%%%%%%%%%%%%%%%%%%%%%%
% Newton iteration with simple continuation
%%%%%%%%%%%%%%%%%%%%%%%%%%%%%%%%%%%%%%%%%%%%%%%%%%%%%%%%%%%%%%%%%%
while ( j_cur <= j & dj > 0 )  

  % display j_cur 
  j_cur = j_cur   

  % initialize Newton iteration
  count     = 0;  % iteration count for single Newton iteration
  c0        = 1-j_cur + epsilon^2*(2*E(1)-2*E(N)-w*(E.^2)); 
  res       = epsilon^2*(L*E-0.5*E.^3) ...
            - 0.25*(c0+j_cur*(x+1)).*E - 0.25*j_cur;
  res(1)    = -k_c*(c0+2*j_cur+epsilon^2*(2*E(1)^2+4*D(1,:)*E)) ...
            + j_r - j_cur;
  res(N)    = k_c*(c0+epsilon^2*(2*E(N)^2+4*D(N,:)*E)) ...
            - j_r - j_cur;
  res_norm  = norm(res,inf);

  while ( (res_norm > res_tol) & (count < max_iters) )
  
    % construct Jacobian for interior grid points
    dc0_dE = -2*epsilon^2*w.*E'; 
    dc0_dE(1) = dc0_dE(1) + 2*epsilon^2;
    dc0_dE(N) = dc0_dE(N) - 2*epsilon^2;
    J = epsilon^2*L ...
      - diag(1.5*epsilon^2*(E.*E) + 0.25*(c0+j_cur*(x+1))) ...
      - 0.25*kron(E,dc0_dE);
    
    % construct Jacobian for boundary conditions
    J(1,:) = -k_c*(dc0_dE + 4*epsilon^2*D(1,:)); 
    J(1,1) = J(1,1) - 4*k_c*epsilon^2*E(1);
    J(N,:) = k_c*(dc0_dE + 4*epsilon^2*D(N,:)); 
    J(N,N) = J(N,N) + 4*k_c*epsilon^2*E(N);
 
    % compute delta_E and update solution
    delta_E = -J\res;  E = E + delta_E;

    % update residual
    c0     = 1-j_cur + epsilon^2*(2*E(1)-2*E(N)-w*(E.^2)); 
    res    = epsilon^2*(L*E-0.5*E.^3) ...
           - 0.25*(c0+j_cur*(x+1)).*E - 0.25*j_cur;
    res(1) = -k_c*(c0+2*j_cur+epsilon^2*(2*E(1)^2+4*D(1,:)*E)) ...
           + j_r - j_cur;
    res(N) = k_c*(c0+epsilon^2*(2*E(N)^2+4*D(N,:)*E)) ...
           - j_r - j_cur;

    % update loop variables
    res_norm  = norm(res,inf)
    count = count + 1
  
  end  % Newton iteration loop

  % update continuation variables
  if (j - j_cur < dj)
    dj = j - j_cur;
  end
  j_cur = j_cur + dj;

end  

%%%%%%%%%%%%%%%%%%%%%%%%%%%%%%%%%%%%%%%%%%%%%%%%%%%%%%%%%%%%%%%%%%
% plot solution
%%%%%%%%%%%%%%%%%%%%%%%%%%%%%%%%%%%%%%%%%%%%%%%%%%%%%%%%%%%%%%%%%%
figure(1); clf;
plot(x,E,'k-');
axis([-1 1 -100 0]);
xlabel('x'); ylabel('E','Rotation',0);

%%%%%%%%%%%%%%%%%%%%%%%%%%%%%%%%%%%%%%%%%%%%%%%%%%%%%%%%%%%%%%%%%%
% plot spectral coefficients
%%%%%%%%%%%%%%%%%%%%%%%%%%%%%%%%%%%%%%%%%%%%%%%%%%%%%%%%%%%%%%%%%%
coefs = abs(fft([E; flipud(E(2:end-1))]));
figure(2); clf;
semilogy(coefs(1:N),'ko');
axis([0 N 1e-15 1e4]);
xlabel('n'); ylabel('|a_n|','Rotation',0,'Position',[-32 5e-6]);
\end{lstlisting}

\newpage
\section{\footnotesize MATLAB Code for Electrochemical Thin-Film Example with 
         Variable Transformation \label{appendix:echem_thin_film_mapped_code}}
\footnotesize This code relies on {\texttt cheb.m} and 
{\texttt clencurt.m}~\cite{trefethen_spectral_book}.
\begin{lstlisting}
%%%%%%%%%%%%%%%%%%%%%%%%%%%%%%%%%%%%%%%%%%%%%%%%%%%%%%%%%%%%%%%%%%
% parameters
%%%%%%%%%%%%%%%%%%%%%%%%%%%%%%%%%%%%%%%%%%%%%%%%%%%%%%%%%%%%%%%%%%
N = 200; beta = 0.75; alpha = atanh(beta);
j = 1.5; epsilon = 0.01; k_c = 10; j_r = 10;    
res_tol = 1e-8; max_iters = 20;

%%%%%%%%%%%%%%%%%%%%%%%%%%%%%%%%%%%%%%%%%%%%%%%%%%%%%%%%%%%%%%%%%%
% compute grid, differentiation matrix, and quadrature weights
%%%%%%%%%%%%%%%%%%%%%%%%%%%%%%%%%%%%%%%%%%%%%%%%%%%%%%%%%%%%%%%%%%
[D,y] = cheb(N-1);                   % Chebyshev differentiation matrix 
[y,w] = clencurt(N-1);               % Clenshaw-Curtis quadrature weights 
x = tanh(alpha*y)/beta;              % mapped grid points
gamma = beta/alpha*cosh(alpha*y).^2; % derivative transformation factor
diag_gamma = diag(gamma);            % cache diag(gamma) 
L = diag_gamma*D*diag_gamma*D;       % transformed Laplacian operator

%%%%%%%%%%%%%%%%%%%%%%%%%%%%%%%%%%%%%%%%%%%%%%%%%%%%%%%%%%%%%%%%%%
% set up continuation in j
%%%%%%%%%%%%%%%%%%%%%%%%%%%%%%%%%%%%%%%%%%%%%%%%%%%%%%%%%%%%%%%%%%
j_start = 0.5; dj = 0.1; j_cur = j_start;  

%%%%%%%%%%%%%%%%%%%%%%%%%%%%%%%%%%%%%%%%%%%%%%%%%%%%%%%%%%%%%%%%%%
% generate initial iterate for Newton iteration
%%%%%%%%%%%%%%%%%%%%%%%%%%%%%%%%%%%%%%%%%%%%%%%%%%%%%%%%%%%%%%%%%%
c0 = 1-j_cur; c = c0 + j_cur*(x+1); E = -2*j_cur./(j_cur*(x+1)+c0);

%%%%%%%%%%%%%%%%%%%%%%%%%%%%%%%%%%%%%%%%%%%%%%%%%%%%%%%%%%%%%%%%%%
% Newton iteration with simple continuation
%%%%%%%%%%%%%%%%%%%%%%%%%%%%%%%%%%%%%%%%%%%%%%%%%%%%%%%%%%%%%%%%%%
while ( j_cur <= j & dj > 0 )  

  % display j_cur 
  j_cur = j_cur   

  % initialize Newton iteration
  count  = 0;  % iteration count for single Newton iteration
  c0     = 1-j_cur ...
         + epsilon^2*(2*gamma(1)*E(1)-2*gamma(N)*E(N)-w*(gamma.*(E.^2))); 
  res    = epsilon^2*(L*(gamma.*E)-0.5*(gamma.^3).*(E.^3)) ...
         - 0.25*(c0+j_cur*(x+1)).*gamma.*E - 0.25*j_cur;
  res(1) = -k_c*(c0+2*j_cur ...
                +epsilon^2*(2*gamma(1)^2*E(1)^2 ...
                           +4*gamma(1)*D(1,:)*(gamma.*E))) ...
         + j_r - j_cur;
  res(N) = k_c*(c0+epsilon^2*(2*gamma(N)^2*E(N)^2 ...
                             +4*gamma(N)*D(N,:)*(gamma.*E))) ...
         - j_r - j_cur;
  res_norm  = norm(res,inf);

  while ( (res_norm > res_tol) & (count < max_iters) )
  
    % construct Jacobian for interior grid points
    dc0_dE = -2*epsilon^2*w.*(gamma.*E)'; 
    dc0_dE(1) = dc0_dE(1) + 2*epsilon^2*gamma(1); 
    dc0_dE(N) = dc0_dE(N) - 2*epsilon^2*gamma(N);
    J = epsilon^2*L*diag_gamma ...
      - diag(1.5*epsilon^2*(gamma.^3).*(E.^2) ...
            +0.25*gamma.*(c0+j_cur*(x+1))) ...
      - 0.25*kron(gamma.*E,dc0_dE);
    
    % construct Jacobian for boundary conditions
    J(1,:) = -k_c*(dc0_dE + 4*epsilon^2*gamma(1)*D(1,:)*diag_gamma); 
    J(1,1) = J(1,1) - 4*k_c*epsilon^2*gamma(1)^2*E(1);
    J(N,:) = k_c*(dc0_dE + 4*epsilon^2*gamma(N)*D(N,:)*diag_gamma); 
    J(N,N) = J(N,N) + 4*k_c*epsilon^2*gamma(N)^2*E(N);
 
    % compute delta_E and update solution
    delta_E = -J\res;  E = E + delta_E;

    % update residual
    c0     = 1-j_cur ...
           + epsilon^2*(2*gamma(1)*E(1)-2*gamma(N)*E(N)-w*(gamma.*(E.^2))); 
    res    = epsilon^2*(L*(gamma.*E)-0.5*(gamma.^3).*(E.^3)) ...
           - 0.25*(c0+j_cur*(x+1)).*gamma.*E - 0.25*j_cur;
    res(1) = -k_c*(c0+2*j_cur ...
                  +epsilon^2*(2*gamma(1)^2*E(1)^2 ...
                             +4*gamma(1)*D(1,:)*(gamma.*E))) ...
           + j_r - j_cur;
    res(N) = k_c*(c0+epsilon^2*(2*gamma(N)^2*E(N)^2 ...
                               +4*gamma(N)*D(N,:)*(gamma.*E))) ...
           - j_r - j_cur;

    % update loop variables
    res_norm  = norm(res,inf)
    count = count + 1
  
  end  % Newton iteration loop

  % update continuation variables
  if (j - j_cur < dj)
    dj = j - j_cur;
  end
  j_cur = j_cur + dj;

end  

%%%%%%%%%%%%%%%%%%%%%%%%%%%%%%%%%%%%%%%%%%%%%%%%%%%%%%%%%%%%%%%%%%
% plot solution
%%%%%%%%%%%%%%%%%%%%%%%%%%%%%%%%%%%%%%%%%%%%%%%%%%%%%%%%%%%%%%%%%%
figure(1); clf;
plot(x,E.*gamma,'k-');
axis([-1 1 -100 0]);
xlabel('x'); ylabel('E','Rotation',0);

%%%%%%%%%%%%%%%%%%%%%%%%%%%%%%%%%%%%%%%%%%%%%%%%%%%%%%%%%%%%%%%%%%
% plot spectral coefficients
%%%%%%%%%%%%%%%%%%%%%%%%%%%%%%%%%%%%%%%%%%%%%%%%%%%%%%%%%%%%%%%%%%
coefs = abs(fft([E; flipud(E(2:end-1))]));
figure(2); clf;
semilogy(coefs(1:N),'ko');
axis([0 N 1e-15 1e4]);
xlabel('n'); ylabel('|a_n|','Rotation',0,'Position',[-25 5e-6]);
\end{lstlisting}

\section{\footnotesize MATLAB Code for Double Layer Charging of Metal Colloid 
         Sphere \label{appendix:metal_colloid_DL_charging_code}}
\footnotesize 
This code relies on {\texttt cheb.m}~\cite{trefethen_spectral_book}.

\begin{lstlisting}
%%%%%%%%%%%%%%%%%%%%%%%%%%%%%%%%%%%%%%%%%%%%%%%%%%%%%%%%%%%%%%%
% Parameters
%%%%%%%%%%%%%%%%%%%%%%%%%%%%%%%%%%%%%%%%%%%%%%%%%%%%%%%%%%%%%%%
% physical parameters
v = 0; E = 10; epsilon = 0.01; delta = 1;

% grid parameters 
N_r = 30;   % number of grid points in radial direction
N_t = 30;   % number of grid points in polar angle direction
L_r = 0.5;  % scale parameter in radial direction

% continuation parameters
E_start = 1; E_final = E; dE = 0.5;

% Newton iteration parameters
res_tol = 1e-8; delta_tol = 1e-13; max_iters = 20;

% zeta-potential iteration parameters
zeta_max_iters = 20; zeta_delta_tol = 1e-13; zeta_res_tol = 1e-9;

% boundary conditions
c_infinity = 1;

%%%%%%%%%%%%%%%%%%%%%%%%%%%%%%%%%%%%%%%%%%%%%%%%%%%%%%%%%%%%%%%
% Construct computational grid and differentiation operators
%%%%%%%%%%%%%%%%%%%%%%%%%%%%%%%%%%%%%%%%%%%%%%%%%%%%%%%%%%%%%%%
% construct the differentiation matrix for the radial coordinate
[D_y,y] = cheb(N_r);
one_minus_y = spdiags(1-y,0,N_r+1,N_r+1);
D_r = 0.5/L_r*(one_minus_y^2)*D_y;
warning off MATLAB:divideByZero
r = L_r*(1+y)./(1-y);
warning on MATLAB:divideByZero
r = r+1; % shift 0 to 1

% construct the differentiation matrix for the polar angle coordinate
theta = (2*[1:N_t]'-1)*pi/2/N_t; T = repmat(theta,1,N_t);
c = ones(1,N_t).*(-1).^([1:N_t]+1);
off_diag_D = ...
  repmat(c,N_t,1).*sin(N_t*T).*sin(T')./(cos(T')-cos(T)+eye(N_t));
diag_D = -0.5*cot(theta);
D_theta = triu(off_diag_D,1) + tril(off_diag_D,-1) + diag(diag_D);

num_gridpts_r = length(r); num_gridpts_theta = length(theta);
num_gridpts = (num_gridpts_r-1)*N_t; 
num_gridpts_interior = num_gridpts-N_t;

% cache common expressions
one_over_r = spdiags(1./r,0,num_gridpts_r,num_gridpts_r);
cos_theta = cos(theta);
sin_theta = sin(theta);
cos_theta_full = kron(cos_theta,ones(num_gridpts_r-1,1));
sin_theta_full = kron(sin_theta,ones(num_gridpts_r-1,1));
sin_theta_mat = spdiags(sin_theta,0,num_gridpts_theta,num_gridpts_theta);
one_over_sin_theta = spdiags(1./sin_theta,0, ...
                             num_gridpts_theta,num_gridpts_theta);

% construct divergence operator
D = {kron(speye(num_gridpts_theta), 2*one_over_r + D_r), ...
     kron(one_over_sin_theta*D_theta*sin_theta_mat, one_over_r)};

% construct gradient operator
G = {kron(speye(num_gridpts_theta),D_r), kron(D_theta,one_over_r)};

% construct laplacian operators
L = kron(speye(num_gridpts_theta),2*one_over_r*D_r + D_r^2) ...
  + kron(one_over_sin_theta*D_theta*sin_theta_mat*D_theta, one_over_r^2);

% construct surface derivative operators
D_s = one_over_sin_theta*D_theta*sin_theta_mat/r(end);
G_s = D_theta/r(end);

% construct normal derivative operator
G_n = -kron(speye(N_t),D_r(end,:));  % d/dn = -d/dr at r = 1

%%%%%%%%%%%%%%%%%%%%%%%%%%%%%%%%%%%%%%%%%%%%%%%%%%%%%%%%%%%%%%%
% Construct matrices to extract subsets of grid points
%%%%%%%%%%%%%%%%%%%%%%%%%%%%%%%%%%%%%%%%%%%%%%%%%%%%%%%%%%%%%%%
% construct matrices to extract the rows corresponding to finite
% grid points (everything except for r = infty)
r_finite_pt_restrictor = spdiags(ones(num_gridpts_r-1,1), 1, ...
                                num_gridpts_r-1, num_gridpts_r);
finite_pt_restrictor = kron(speye(N_t),r_finite_pt_restrictor);

% construct matrix to extract the rows corresponding to interior
% grid points (everything except for r = 1 and r = infty)
r_interior_restrictor = spdiags(ones(num_gridpts_r-2,1), 1, ...
                               num_gridpts_r-2, num_gridpts_r);
interior_restrictor = kron(speye(N_t),r_interior_restrictor);

% construct matrix to extract the rows corresponding to r = 1 (surface)
% from a vector that already has r = infinity removed
r_surf_restrictor = spalloc(1,num_gridpts_r-1,1);
r_surf_restrictor(1,num_gridpts_r-1) = 1;
surf_restrictor = kron(speye(N_t),r_surf_restrictor);

% construct matrix to extract the rows corresponding to r = infty
r_inf_restrictor = spalloc(1,num_gridpts_r,1);
r_inf_restrictor(1,1) = 1;
inf_restrictor = kron(speye(N_t),r_inf_restrictor);

% extract part of G operator that contributes to finite points 
% using finite points
G_f = {finite_pt_restrictor*G{1}*finite_pt_restrictor', ...
       finite_pt_restrictor*G{2}*finite_pt_restrictor'};

% split the D operators into two parts:
% (1) contributions from finite points to finite points
% (2) contributions from infinity to finite points
D_f = {interior_restrictor*D{1}*finite_pt_restrictor', ...
       interior_restrictor*D{2}*finite_pt_restrictor'};
D_inf = {interior_restrictor*D{1}*inf_restrictor', ...
         interior_restrictor*D{2}*inf_restrictor'};

% split the G_n operators into two parts:
% (1) contributions from finite points
% (2) contributions from infinity
G_n_f = G_n*finite_pt_restrictor';
G_n_inf = G_n*inf_restrictor';

% split the Laplacian operators into two parts:
% (1) contributions from finite points to finite points
% (2) contributions from infinity to finite points
L_f = interior_restrictor*L*finite_pt_restrictor';
L_inf = interior_restrictor*L*inf_restrictor';

%%%%%%%%%%%%%%%%%%%%%%%%%%%%%%%%%%%%%%%%%%%%%%%%%%%%%%%%%%%%%%%
% Continuation loop for Newton iteration
%%%%%%%%%%%%%%%%%%%%%%%%%%%%%%%%%%%%%%%%%%%%%%%%%%%%%%%%%%%%%%%
E = E_start; c = ones(num_gridpts,1);  psi = zeros(num_gridpts,1);  
while ( E <= E_final & dE > 0 )

  % Show progress information
  mesg = sprintf('E = %f', E); disp(mesg);

  %%%%%%%%%%%%%%%%%%%%%%%%%%%%%%%%%%%%%%%%%%%%%%%%%
  % compute constant terms in F = (F1,F2,H1,H2)
  %%%%%%%%%%%%%%%%%%%%%%%%%%%%%%%%%%%%%%%%%%%%%%%%%
  F1_const_term = c_infinity*(L_inf*ones(N_t,1));
  F2_const_term = E*c_infinity*(-D_inf{1}*cos_theta + D_inf{2}*sin_theta);
  H2_const_term = - c_infinity*(G_n_inf*ones(N_t,1));
  
  %%%%%%%%%%%%%%%%%%%%%%%%%%%%%%%%%%%%%%%%%%%%%%%%%
  % compute constant parts of Jacobian
  %%%%%%%%%%%%%%%%%%%%%%%%%%%%%%%%%%%%%%%%%%%%%%%%%
  DF1_Dc_const = L_f;
  DF2_Dc_const = ...
    - D_f{1}*spdiags(E*cos_theta_full,0,num_gridpts,num_gridpts) ...
    + D_f{2}*spdiags(E*sin_theta_full,0,num_gridpts,num_gridpts);
  DH1_Dc_const = spalloc(N_t,num_gridpts,N_t);
  DH1_Dc_const(:,num_gridpts_r-1:num_gridpts_r-1:end) = ...
    -spdiags(E*cos_theta,0,N_t,N_t);
  DH2_Dc_const = -G_n_f;
  
  %%%%%%%%%%%%%%%%%%%%%%%%%%%%%%%%%%%%%%%%%%%%%%%%%
  % initialize loop variables using current 
  % solution for c and psi
  %%%%%%%%%%%%%%%%%%%%%%%%%%%%%%%%%%%%%%%%%%%%%%%%%
  % extract surface concentration and potential
  c_s = c(num_gridpts_r-1:num_gridpts_r-1:end);
  phi_s = psi(num_gridpts_r-1:num_gridpts_r-1:end) - E*cos_theta;
  
  % compute zeta potential
  zeta = computeZetaPotential( ...
    v-phi_s, c_s, delta, zeta_res_tol, zeta_delta_tol, zeta_max_iters); 
  
  % cache some common expressions
  log_c_s = log(c_s);
  sinh_zeta_over_two = sinh(zeta/2);
  cosh_zeta_over_two = cosh(zeta/2);

  % compute surface charge density and excess neutral ion concentration
  q = -2*sqrt(c_s).*sinh_zeta_over_two;
  w = 4*sqrt(c_s).*(sinh(zeta/4)).^2;
  
  % compute initial residual
  F1 = F1_const_term + L_f*c;
  F2 = F2_const_term + D_f{1}*(c.*(G_f{1}*psi-E*cos_theta_full)) ...
                     + D_f{2}*(c.*(G_f{2}*psi+E*sin_theta_full)); 
  H1 = epsilon*D_s*(q.*(G_s*log_c_s) + w.*(G_s*phi_s)) ...
     - c_s.*(G_n_f*psi + E*cos_theta);
  H2 = epsilon*D_s*(w.*(G_s*log_c_s) + q.*(G_s*phi_s)) ...
     - G_n_f*c + H2_const_term;
  F = [F1; F2; H1; H2];  
  res = norm(F,inf);
  
  %%%%%%%%%%%%%%%%%%%%%%%%%%%%%%%%%%%%%%%%%%%%%%%%%
  % Newton iteration loop
  %%%%%%%%%%%%%%%%%%%%%%%%%%%%%%%%%%%%%%%%%%%%%%%%%
  norm_delta_soln = 1;
  count = 0;

  % begin Newton iteration loop
  while (res > res_tol && norm_delta_soln > delta_tol && count < max_iters)
  
    % compute Jacobian
    Dzeta_Dpsi = -1./(1+delta*sqrt(c_s).*cosh_zeta_over_two);
    Dzeta_Dc_s = -delta*sinh_zeta_over_two./sqrt(c_s) ...
               ./(1+delta*sqrt(c_s).*cosh_zeta_over_two);
    DH1_Dc_var = ( epsilon * D_s * ( ...
       + spdiags(0.5*q./c_s.*(G_s*log_c_s) ...
                -sqrt(c_s).*cosh_zeta_over_two.*Dzeta_Dc_s.*(G_s*log_c_s), ...
                 0,N_t,N_t) ...
       + spdiags(q,0,N_t,N_t)*G_s*spdiags(1./c_s,0,N_t,N_t) ...
       + spdiags(0.5*w./c_s.*(G_s*phi_s) ...
                +sqrt(c_s).*sinh_zeta_over_two.*Dzeta_Dc_s.*(G_s*phi_s), ...
                 0,N_t,N_t) ) ...
       - spdiags(G_n_f*psi,0,N_t,N_t) )*surf_restrictor;
    DH1_Dpsi_var = epsilon * D_s * ( ...
       - spdiags(sqrt(c_s).*cosh_zeta_over_two.*Dzeta_Dpsi.*(G_s*log_c_s), ...
                 0,N_t,N_t) ...
       + spdiags(w,0,N_t,N_t)*G_s ...
       + spdiags(sqrt(c_s).*sinh_zeta_over_two.*Dzeta_Dpsi.*(G_s*phi_s), ...
                 0,N_t,N_t) ) * surf_restrictor ...
       - spdiags(c_s,0,N_t,N_t)*G_n_f;
    DH2_Dc_var = epsilon * D_s * ( ...
       + spdiags(0.5*w./c_s.*(G_s*log_c_s) ...
                +sqrt(c_s).*sinh_zeta_over_two.*Dzeta_Dc_s.*(G_s*log_c_s), ...
                 0,N_t,N_t) ...
       + spdiags(w,0,N_t,N_t)*G_s*spdiags(1./c_s,0,N_t,N_t) ...
       + spdiags(0.5*q./c_s.*(G_s*phi_s) ...
                -sqrt(c_s).*cosh_zeta_over_two.*Dzeta_Dc_s.*(G_s*phi_s), ...
                 0,N_t,N_t) ) * surf_restrictor;
    DH2_Dpsi_var = epsilon * D_s * ( ...
         spdiags(sqrt(c_s).*sinh_zeta_over_two.*Dzeta_Dpsi.*(G_s*log_c_s), ...
                 0,N_t,N_t) ...
       + spdiags(q,0,N_t,N_t)*G_s ...
       - spdiags(sqrt(c_s).*cosh_zeta_over_two.*Dzeta_Dpsi.*(G_s*phi_s), ...
                 0,N_t,N_t) ) * surf_restrictor; 
    J = [DF1_Dc_const, spalloc(num_gridpts_interior,num_gridpts,0); ...
         ( DF2_Dc_const ...
         + D_f{1}*spdiags(G_f{1}*psi,0,num_gridpts,num_gridpts) ...
         + D_f{2}*spdiags(G_f{2}*psi,0,num_gridpts,num_gridpts) ), ...
         ( D_f{1}*spdiags(c,0,num_gridpts,num_gridpts)*G_f{1} ...
         + D_f{2}*spdiags(c,0,num_gridpts,num_gridpts)*G_f{2} ); ...
         ( DH1_Dc_const + DH1_Dc_var ), DH1_Dpsi_var; ...
         ( DH2_Dc_const + DH2_Dc_var ), DH2_Dpsi_var]; 

    % compute delta_soln
    delta_soln = -J\F;
  
    % update solution 
    c = c + delta_soln(1:num_gridpts);
    psi = psi + delta_soln(num_gridpts+1:end);
  
    %%%%%%%%%%%%%%%%%%%%%%%%%%%%
    % update residual
    %%%%%%%%%%%%%%%%%%%%%%%%%%%%
    % extract surface concentration and potential
    c_s = c(num_gridpts_r-1:num_gridpts_r-1:end);
    phi_s = psi(num_gridpts_r-1:num_gridpts_r-1:end) - E*cos_theta;
  
    % compute zeta potential
    zeta = computeZetaPotential( ...
      v-phi_s, c_s, delta, zeta_res_tol, zeta_delta_tol, zeta_max_iters); 
  
    % cache some common expressions
    log_c_s = log(c_s);
    sinh_zeta_over_two = sinh(zeta/2);
    cosh_zeta_over_two = cosh(zeta/2);

    % compute surface charge density and excess neutral ion concentration
    q = -2*sqrt(c_s).*sinh_zeta_over_two;
    w = 4*sqrt(c_s).*(sinh(zeta/4)).^2;
  
    % compute residual
    F1 = F1_const_term + L_f*c;
    F2 = F2_const_term + D_f{1}*(c.*(G_f{1}*psi-E*cos_theta_full)) ...
                       + D_f{2}*(c.*(G_f{2}*psi+E*sin_theta_full)); 
    H1 = epsilon*D_s*(q.*(G_s*log_c_s) + w.*(G_s*phi_s)) ...
       - c_s.*(G_n_f*psi + E*cos_theta);
    H2 = epsilon*D_s*(w.*(G_s*log_c_s) + q.*(G_s*phi_s)) ...
       - G_n_f*c + H2_const_term;
    F = [F1; F2; H1; H2];  
    res = norm(F,inf);

    % update norm_delta_soln, count, and residual history
    norm_delta_soln = norm(delta_soln,inf);
    count = count + 1;
  
    % show stats
    status = [res norm_delta_soln count]
  
  end  % end Newton iteration loop

  % update E
  if (E_final - E< dE)
    dE = E_final - E;
  end
  E = E + dE;

end  
  
%%%%%%%%%%%%%%%%%%%%%%%%%%%%%%%%%%%%%%%%%%%%%%%%%%%%%%%%%%%%%%%
% Append values at infinity to results
%%%%%%%%%%%%%%%%%%%%%%%%%%%%%%%%%%%%%%%%%%%%%%%%%%%%%%%%%%%%%%%
c = finite_pt_restrictor'*c; 
c(1:num_gridpts_r:end) = c_infinity;
psi = finite_pt_restrictor'*psi; 

%%%%%%%%%%%%%%%%%%%%%%%%%%%%%%%%%%%%%%%%%%%%%%%%%%%%%%%%%%%%%%%%%%
%  Plot results
%%%%%%%%%%%%%%%%%%%%%%%%%%%%%%%%%%%%%%%%%%%%%%%%%%%%%%%%%%%%%%%%%%
axis_scale = 3;

% psi = potential relative to applied field
figure(1); clf;
min_psi = min(psi); max_psi = max(psi);
x_scale = axis_scale; y_scale = axis_scale/2;
[rr,tt] = meshgrid(r,pi/2-theta); [xx,yy] = pol2cart(tt,rr);
surf(xx,yy,reshape(psi',N_r+1,N_t)');
hold on;
surf_theta = pi/2-[0; theta; pi];
[surf_x,surf_y] = pol2cart(surf_theta,ones(size(surf_theta)));
surf_z = min_psi*ones(size(surf_x));
proj_color = [180 200 220]/256;
plot3(surf_x,surf_y,surf_z,'k'); fill3(surf_x,surf_y,surf_z,proj_color);
xlabel('x'); ylabel('z'); zlabel('\psi','rotation',0);
axis([0 x_scale -y_scale y_scale min_psi max_psi]);

% concentration
figure(2); clf;
min_c = min(c); max_c = max(c);
x_scale = axis_scale; y_scale = axis_scale/2;
[rr,tt] = meshgrid(r,pi/2-theta); [xx,yy] = pol2cart(tt,rr);
surf(xx,yy,reshape(c',N_r+1,N_t)');
hold on;
surf_theta = pi/2-[0; theta; pi];
[surf_x,surf_y] = pol2cart(surf_theta,ones(size(surf_theta)));
surf_z = min_c*ones(size(surf_x));
proj_color = [180 200 220]/256;
plot3(surf_x,surf_y,surf_z,'k'); fill3(surf_x,surf_y,surf_z,proj_color);
xlabel('x'); ylabel('z'); zlabel('c','rotation',0,'position',[-1.4 2 0.85]);
axis([0 x_scale -y_scale y_scale min_c max_c]);
\end{lstlisting}

\subsection{\footnotesize \texttt{computeZetaPotential()}}
\begin{lstlisting}
function zeta = computeZetaPotential(...
  Psi, c_s, delta, res_tol, delta_zeta_tol, max_iters)

% initialize iteration 
zeta = Psi;  % use Psi as an initial guess for zeta
delta_zeta = 1; norm_delta_zeta = norm(delta_zeta,inf);
res = 1; norm_res = norm(res,inf);
count = 0;
res = zeta + 2*delta*sqrt(c_s).*sinh(zeta/2) - Psi;

% Newton iteration
while (norm_res > res_tol & norm_delta_zeta > delta_zeta_tol ...
       & count < max_iters)

  J = 1 + delta*sqrt(c_s).*cosh(zeta/2);
  delta_zeta = -res./J;
  zeta = zeta + delta_zeta;
  res = zeta + 2*delta*sqrt(c_s).*sinh(zeta/2) - Psi;
  norm_res = norm(res,inf);
  norm_delta_zeta= norm(delta_zeta,inf);
  count = count + 1;

end
\end{lstlisting}

\end{multicols}
\end{landscape}

\end{document}